\documentclass[11pt,a4paper]{article}

\usepackage[T1]{fontenc}
\usepackage[utf8]{inputenc}
\usepackage[english]{babel}
\usepackage{lmodern}

\usepackage{amsmath,amssymb,amsthm,mathtools}

\usepackage[margin=1in]{geometry}
\usepackage{microtype}
\usepackage{booktabs}

\usepackage[hidelinks]{hyperref}
\newtheorem{theorem}{Theorem}[section]
\newtheorem{proposition}[theorem]{Proposition}
\newtheorem{lemma}[theorem]{Lemma}
\newtheorem{corollary}[theorem]{Corollary}

\theoremstyle{definition}
\newtheorem{definition}[theorem]{Definition}

\theoremstyle{remark}
\newtheorem{rmk}[theorem]{Remark}

\newenvironment{pf}[1][\proofname]
  {\begin{proof}[#1]}
  {\end{proof}}

\title{Sharp Simultaneous Asymptotic Expansions for Sample-Mean
Operators from Natural Exponential Families}

\author{%
  Oleksandr Volkov\\[0.4em]
  \small Department of Statistics, University of California at Berkeley  
}

\date{}

\makeatletter
\renewcommand{\@seccntformat}[1]{\csname the#1\endcsname.\quad}
\makeatother

\usepackage{fancyhdr}
\newcommand{\shorttitle}[1]{\fancyhead[C]{\small #1}}
\shorttitle{Sharp Expansions for Sample-Mean Operators}

\begin{document}

\maketitle

\begin{abstract}
We study the asymptotic behavior of positive linear operators
constructed from sample means of natural exponential families
parametrized by their means. We show that the variance function
\(V\) determines the entire hierarchy of cumulants, central
moments, and coefficient differential operators
\(\mathcal A_r[V]\) governing the asymptotic expansions.
The cumulants satisfy the recursion
\(\kappa_{j+1}=V\kappa_j'\), whereas the central moments are
expressed in terms of partial Bell polynomials. Quantitative
expansions of arbitrary order, uniform on compact subsets of
the mean domain, are obtained. We also establish simultaneous
expansions for the derivatives \((L_nf)^{(\ell)}\), with an
explicit remainder estimate in terms of the modulus of
continuity of the highest derivative. We prove that
\(C^{2N+\ell}\)-smoothness is sufficient for the expansion
of the derivative of order \(\ell\) up to and including the
term of order \(n^{-N}\). The Gauss--Weierstrass operators
show that, within the scale of spaces \(C^m\), this smoothness
threshold cannot generally be lowered and that the order of
the modulus-of-continuity term is optimal. The general
framework is applied to the Bernstein, Sz\'asz--Mirakjan,
Baskakov, Post--Widder, and Gauss--Weierstrass operators.
\end{abstract}

\par


\section{Introduction}
\label{sec:introduction}

Many classical positive approximation operators admit a
probabilistic interpretation as averaging operators. In its
simplest form,
\(
(L_nf)(x)=\mathbb E_x f(\eta_{n,x}),
\)
where, as \(n\to\infty\), the random variable \(\eta_{n,x}\)
concentrates near the point \(x\). This scheme includes the
Bernstein, Sz\'asz--Mirakjan, Baskakov, Post--Widder, and
Gauss--Weierstrass operators. A common structural feature of
these examples is the differential identity
\[
\frac{V(x)}{n}(L_nf)'(x)
=
\mathbb E_x\bigl[(\eta_{n,x}-x)f(\eta_{n,x})\bigr],
\]
where \(V\) is a positive function. Operators satisfying
identities of this type were studied by May, later by Ismail
and May
\cite{May1976,IsmailMay1978}; such operators are
usually referred to as operators of exponential type.

The same structure arises naturally in probability theory
and mathematical statistics. Let
\(\{Q_x:x\in\mathcal M\}\) be a natural exponential family
parametrized by its mean, and let
\(
V(x)=\operatorname{Var}_x(\xi_x),
\) 
\(
\xi_x\sim Q_x,
\)
be its variance function. Variance functions play a
fundamental role in the theory of natural exponential
families; in particular, families with quadratic variance
functions were classified by Morris~\cite{Morris1982}. If
\(\eta_{n,x}\) is the arithmetic mean of \(n\) independent
copies of the random variable \(\xi_x\), then the family of
its distributions
\(\{Q_{n,x}:x\in\mathcal M\}\) again forms a natural
exponential family --- in general, with a different base
measure --- with mean \(x\) and variance function \(V(x)/n\).
The corresponding sample-mean operator
\((L_nf)(x)=\mathbb E_x f(\eta_{n,x})\) satisfies the defining
differential identity of operators of exponential type. In
this sense, sample-mean operators provide a natural
probabilistic realization of the differential structure of
exponential type and unify the classical approximation
operators listed above within a single scheme, without
requiring the function \(V\) to be quadratic.

The asymptotic study of positive linear operators goes back
to the classical Voronovskaya formula
\cite{Voronovskaya1932}. General asymptotic expansions in
terms of central moments were developed by
Sikkema~\cite{Sikkema1970}, whereas the differential
structure and recurrence relations for the moments of
operators of exponential type were studied in
\cite{IsmailMay1978}.
Simultaneous approximation, in which the derivatives of
\(L_nf\) are approximated together with the function \(f\),
has been investigated both for individual classical
operators and in more general settings. For Bernstein
operators, Floater~\cite{Floater2005} obtained an error
estimate and a differentiated first-order asymptotic formula
under the assumption \(f\in C^{\ell+2}[0,1]\) for the
derivative of order \(\ell\). A related problem concerning
the compatibility of differentiation with passage to the
limit in Voronovskaya-type expansions for Bernstein-type
operators was studied by Xiang~\cite{Xiang2022}. General
results on simultaneous approximation for positive linear
operators and operators of exponential type were obtained,
in particular, by Gonska and P\u{a}lt\u{a}nea
\cite{GonskaPaltanea2010} and by
Holho\c{s}~\cite{Holhos2019,Holhos2022Exponential,Holhos2022Simultaneous}. More recently,
Abel~\cite{Abel2023} proved that complete asymptotic
expansions can be differentiated term by term for a class
containing operators of exponential type. In the case of
operators of exponential type, his pointwise expansion of
the \(\ell\)-th derivative through the term of order
\(n^{-N}\) is established under a smoothness assumption of
order \(2(N+\ell)\).

In contrast to general results on termwise differentiation
of asymptotic expansions, we exploit the special structure
of sample means of natural exponential families to construct
explicit coefficient operators, prove a
quantitative estimate uniform on compact subsets with the
smoothness threshold \(C^{2N+\ell}\), and establish the
sharpness of this threshold in the scale of integer-order
smoothness classes \(C^m\). The key point is that the
variance function determines not only the concentration rate
of \(\eta_{n,x}\), but also the entire system of expansion
coefficients. The cumulants of the family satisfy the
recursion
\(\kappa_1=x\),
\(\kappa_{j+1}=V\kappa_j'\), \(j\geq1\), while exact formulas
for the central moments are expressed in terms of partial
Bell polynomials. This leads to differential operators
\(\mathcal A_r[V]\) whose coefficients are differential
polynomials in \(V\). For each \(r\), the operator
\(\mathcal A_r[V]\) has exact order \(2r\), and its leading
coefficient is \(V^r/(2^r r!)\).

The main result is a quantitative expansion for the
derivatives of the operators \(L_n\), uniform on compact
subsets. Let
\(\ell,m\in\mathbb N_0\), \(m\geq\ell\), and
\(s_\ell=\lfloor(m-\ell)/2\rfloor\). For each compact set
\(J\subset\mathcal M\) and its fixed compact neighborhood
\(J_\rho\Subset\mathcal M\), we prove the estimate
\[
\begin{aligned}
\sup_{x\in J}\left|
(L_nf)^{(\ell)}(x)
-
\sum_{r=0}^{s_\ell}
\frac{1}{n^r}
\bigl(\mathcal A_r[V]f\bigr)^{(\ell)}(x)
\right|
\leq{}&
C n^{-(m-\ell)/2}
\omega_{m,J_\rho}\bigl(f;n^{-1/2}\bigr)
+D n^{-s_\ell-1}.
\end{aligned}
\]
Here \(\omega_{m,J_\rho}(f;\delta)\) denotes the modulus of
continuity of the derivative \(f^{(m)}\) on \(J_\rho\), and
the constants \(C,D\) are independent of \(n\). The
smoothness of \(f\) is required only in a neighborhood of the
compact set, whereas outside this neighborhood Borel
measurability and at most polynomial growth suffice. In
particular, the condition
\(f\in C^{2N+\ell}(\mathcal M)\) ensures a locally uniform
expansion of \((L_nf)^{(\ell)}\) through the term of order
\(n^{-N}\), with remainder \(o(n^{-N})\). Thus, each
successive order of the expansion requires two additional
derivatives, whereas each simultaneous differentiation
requires only one. Compared with the condition
\(C^{2(N+\ell)}\) in the general result of
Abel~\cite{Abel2023}, the score-polynomial representation
\eqref{eq:score-representation}, combined with the exact
polynomial identity
\eqref{eq:exact-polynomial-expansion}, yields the sufficient
smoothness threshold \(C^{2N+\ell}\), thereby saving
\(\ell\) derivatives.

This smoothness threshold is sharp in the scale of the
spaces \(C^m\). The principal part of the operator
\(\bigl(\mathcal A_N[V]f\bigr)^{(\ell)}\) contains the
nonzero term
\(V^Nf^{(2N+\ell)}/(2^N N!)\). For the Gauss--Weierstrass
operators, we construct compactly supported functions
\(f\in C^{2N+\ell-1}(\mathbb R)\) for which
\[
n^N\left|
(W_nf)^{(\ell)}(0)
-
\sum_{r=0}^{N-1}\frac{1}{n^r}
\bigl(\mathcal A_r[1/2]f\bigr)^{(\ell)}(0)
\right|
\longrightarrow\infty.
\]
Consequently, even after subtracting all coefficient terms
defined on the class \(C^{2N+\ell-1}\), the remainder is not
of order \(o(n^{-N})\). Another family of Gaussian examples
shows that the term involving the modulus of continuity in
the quantitative estimate above is also sharp in order.
Thus, the paper establishes not only asymptotic expansions
of arbitrary order, but also the sharpness of the smoothness
threshold and of the order of the modulus-of-continuity term
in the quantitative simultaneous estimate.

The proof relies on a representation of the derivatives
\((L_nf)^{(\ell)}\) in terms of logarithmic-derivative
polynomials of the distribution and on an exact identity for
the action of \(L_n\) on polynomials. This makes it possible
to identify the derivatives of the coefficients without
differentiating the asymptotic remainder. The general scheme
covers all five classical families listed above and yields
explicit ordinary and simultaneous expansions for them.

In Section~\ref{sec:exponential_families}, we construct the
exponential representation and introduce the sample-mean
operators. In Section~\ref{sec:cumulants-and-moments}, we
derive formulas for the cumulants, central moments, and
coefficient operators. Sections
\ref{sec:all-order-expansion}--\ref{sec:simultaneous-consequences}
contain the ordinary and simultaneous expansions and their
consequences, while Section~\ref{sec:sharpness} establishes
the sharpness of the smoothness threshold and of the
quantitative estimate. Classical examples are considered in
Section~\ref{sec:classical-operators}.

\section{Natural Exponential Families and Their Sample-Mean Operators}
\label{sec:exponential_families}

\subsection{Mean Parametrization and the Variance Function}
\label{subsec_mean_parametrization}

Let \(\mu\) be a nonzero positive Borel measure on
\(\mathbb R\), not concentrated at a single point, and suppose
that
\[
\Theta
=
\operatorname{int}
\left\{
\theta\in\mathbb R:
\int_{\mathbb R}e^{\theta t}\,\mu(dt)<\infty
\right\}
\]
is a nonempty open interval. Set
\[
K(\theta)
=
\log\int_{\mathbb R}e^{\theta t}\,\mu(dt),
\qquad \theta\in\Theta.
\]
Then
\[
P_\theta(dt)
=
\exp\{\theta t-K(\theta)\}\,\mu(dt),
\qquad \theta\in\Theta,
\]
forms a natural exponential family. If
\(X\sim P_\theta\), then \(K\) is real-analytic on
\(\Theta\), and
\[
K'(\theta)=\mathbb E_\theta X,
\qquad
K''(\theta)=\operatorname{Var}_\theta(X)>0.
\]
Thus, \(K'\) is strictly increasing. Denote by
\(
\mathcal M=K'(\Theta)
\)
the mean domain, and let \(\theta=\theta(x)\) be the inverse
function of \(K'\). Under the mean parametrization, the
family takes the form
\begin{equation}
\label{eq_mean_parametrized_nef}
Q_x(dt)
=
\exp\{\theta(x)t-K(\theta(x))\}\,\mu(dt),
\qquad x\in\mathcal M.
\end{equation}
If \(\xi_x\sim Q_x\), then
\(
\mathbb E_x\xi_x=x,
\)
and the variance function of the family is given by
\[
V(x)
=
\operatorname{Var}_x(\xi_x)
=
K''(\theta(x)),
\qquad x\in\mathcal M.
\]
The identity \(K'(\theta(x))=x\) implies that
\(
\theta'(x)=\frac{1}{V(x)}.
\)

For general background on natural exponential families,
their mean parametrization, and variance functions, see
\cite{BarndorffNielsen1978,Letac1992}.

Fix \(x_0\in\mathcal M\) and set
\(
\tau(x)
=
\int_{x_0}^{x}\frac{du}{V(u)},
\) 
\(
H(x)
=
\int_{x_0}^{x}\frac{u\,du}{V(u)}.
\)
Then
\begin{equation}
\label{eq_tau_H_relations}
\tau(x)=\theta(x)-\theta(x_0),
\qquad
H(x)=K(\theta(x))-K(\theta(x_0)).
\end{equation}

\subsection{Variance-Function Representation and Differential Identities}
\label{subsec_differential_identities}

Let \(\nu=Q_{x_0}\). From
\eqref{eq_mean_parametrized_nef} and
\eqref{eq_tau_H_relations}, we obtain the relative
exponential representation
\begin{equation}
\label{eq:relative_exponential_representation}
Q_x(dt)
=
\exp\{\tau(x)t-H(x)\}\,\nu(dt),
\qquad x\in\mathcal M.
\end{equation}

Denote by
\(
M(z,x)
=
\mathbb E_xe^{z\xi_x}
\)
the moment generating function of the distribution \(Q_x\),
and set
\[
\mathcal D
=
\bigl\{
(z,x)\in\mathbb R\times\mathcal M:
\tau(x)+z\in\tau(\mathcal M)
\bigr\}.
\]
For \((z,x)\in\mathcal D\),
\begin{equation}
\label{eq:mgf_representation}
M(z,x)
=
\exp\left\{
H\left(\tau^{-1}(\tau(x)+z)\right)-H(x)
\right\}.
\end{equation}
Indeed, this formula is equivalent to the identity
\[
M(z,x)
=
\exp\{K(\theta(x)+z)-K(\theta(x))\}.
\]
Differentiating \eqref{eq:mgf_representation} gives
\begin{equation}
\label{eq:differential_characterization}
M_z(z,x)-V(x)M_x(z,x)
=
xM(z,x),
\qquad (z,x)\in\mathcal D.
\end{equation}

The following identity is an integral form of this
differential characterization.

\begin{corollary}
\label{cor:weak_differential_identity}
Let \(f\colon\mathbb R\to\mathbb C\) be a Borel measurable function of
at most polynomial growth. Then
\begin{equation}
\label{eq:weak_differential_identity}
V(x)\frac{d}{dx}
\int_{\mathbb R}f(t)\,Q_x(dt)
=
\int_{\mathbb R}(t-x)f(t)\,Q_x(dt),
\qquad x\in\mathcal M.
\end{equation}
\end{corollary}

\begin{pf}
The finiteness of the Laplace transform of \(\nu\) in a
neighborhood of the point \(\tau(x)\) allows us to
differentiate \eqref{eq:relative_exponential_representation}
under the integral sign. Since
\(
\tau'(x)=\frac{1}{V(x)},
\) 
\(
H'(x)=\frac{x}{V(x)},
\)
we obtain \eqref{eq:weak_differential_identity}.
\end{pf}

To clarify the role of
\eqref{eq:differential_characterization}, we give a concise
converse statement.

\begin{proposition}
\label{prop:differential_characterization}
Let \(\mathcal M\subset\mathbb R\) be a nonempty open
interval, let \(V\colon\mathcal M\to(0,\infty)\) be
continuous, let \(x_0\in\mathcal M\), and let \(\tau,H\) and
\(\mathcal D\) be defined as above.
Suppose that the family of probability measures
\(\{Q_x:x\in\mathcal M\}\) has a finite moment generating
function \(M\in C^1(\mathcal D)\) satisfying
\eqref{eq:differential_characterization} and the condition
\(
M(0,x)=1,
\) 
\( x\in\mathcal M.
\)
Then, for some probability measure \(\nu\), the family admits
the representation
\eqref{eq:relative_exponential_representation}, and \(M\) is
given by \eqref{eq:mgf_representation}. In particular, the
family is a restriction of a natural exponential family
parametrized by the mean \(x\), with variance function \(V\).
\end{proposition}

\begin{pf}
The characteristic system for
\eqref{eq:differential_characterization} is
\(
\frac{dz}{dr}=1,
\)
\(
\frac{dx}{dr}=-V(x),
\)
and along the characteristics the quantities
\(
\tau(x)+z,
\)
\(
e^{H(x)}M(z,x)
\)
are constant. The characteristic through
\((z,x)\in\mathcal D\) intersects the line \(z=0\) at
\(
\left(0,\tau^{-1}(\tau(x)+z)\right),
\)
and therefore the condition \(M(0,x)=1\) gives
\eqref{eq:mgf_representation}.

Set \(\nu=Q_{x_0}\). Then, for
\(s\in\tau(\mathcal M)\),
\[
\int_{\mathbb R}e^{st}\,\nu(dt)
=
M(s,x_0)
=
\exp\left\{H\left(\tau^{-1}(s)\right)\right\}.
\]
The measure defined by the right-hand side of
\eqref{eq:relative_exponential_representation} has the same
moment generating function as \(Q_x\) in a neighborhood of
zero and therefore coincides with \(Q_x\). Finally, setting
\(
\widetilde K(s)=H(\tau^{-1}(s)),
\)
we obtain
\[
\widetilde K'(\tau(x))=x,
\qquad
\widetilde K''(\tau(x))=V(x),
\]
which proves the final assertion. 
\end{pf}

\begin{rmk}
\label{rem:maximal_natural_parameter_domain}
If \(\Theta_\nu\) is the maximal open domain on which the
Laplace transform of \(\nu\) is finite, then we obtain only
\(\tau(\mathcal M)\subseteq\Theta_\nu\). The constructed
family is full if and only if
\(\tau(\mathcal M)=\Theta_\nu\). The question of which
positive functions can occur as variance functions of
natural exponential families is not considered here; for a
survey of this problem and, in particular, of families with
quadratic variance functions, see
\cite{Letac1992,Morris1982}.
\end{rmk}

\subsection{Sample Means and Associated Positive Operators}
\label{subsec_sample_means}

For \(x\in\mathcal M\), let
\(\xi_{1,x},\xi_{2,x},\ldots\) be independent identically
distributed random variables with distribution \(Q_x\), and
let
\(
\eta_{n,x}
=
\frac{1}{n}\sum_{j=1}^{n}\xi_{j,x}.
\)
Denote the distribution of \(\eta_{n,x}\) by \(Q_{n,x}\).
Then
\begin{equation}
\label{eq:sample_mean_expectation_variance}
\mathbb E_x\eta_{n,x}=x,
\qquad
\operatorname{Var}_x(\eta_{n,x})
=
\frac{V(x)}{n}.
\end{equation}

Let \(\nu_n\) be the pushforward of the measure
\(\nu^{\otimes n}\) under the mapping
\[
T_n(t_1,\ldots,t_n)
=
\frac{t_1+\cdots+t_n}{n}.
\]

\begin{proposition}
\label{prop:sample_mean_family}
For each \(n\geq1\), the family
\(\{Q_{n,x}:x\in\mathcal M\}\) is a restriction of a natural
exponential family parametrized by the mean \(x\), with
variance function \(V(x)/n\). Moreover,
\begin{equation}
\label{eq:sample_mean_exponential_representation}
Q_{n,x}(dt)
=
\exp\{n\tau(x)t-nH(x)\}\,\nu_n(dt),
\end{equation}
and its moment generating function is given by
\begin{equation}
\label{eq:sample_mean_mgf}
\begin{aligned}
M_n(z,x)
&=
\mathbb E_xe^{z\eta_{n,x}}
=
\left(M\left(\frac{z}{n},x\right)\right)^n
=
\exp\left\{
n\left[
H\left(\tau^{-1}\left(\tau(x)+\frac{z}{n}\right)\right)
-H(x)
\right]
\right\},
\end{aligned}
\end{equation}
where \(\tau(x)+z/n\in\tau(\mathcal M)\).
\end{proposition}

\begin{pf}
It follows from
\eqref{eq:relative_exponential_representation} that
\[
\prod_{j=1}^{n}Q_x(dt_j)
=
\exp\left\{
n\tau(x)T_n(t_1,\ldots,t_n)-nH(x)
\right\}
\nu^{\otimes n}(dt_1\cdots dt_n).
\]
Taking the pushforward under \(T_n\) gives
\eqref{eq:sample_mean_exponential_representation}. Together
with \eqref{eq:sample_mean_expectation_variance}, this proves
the assertion concerning the parametrization and the
variance function. Formula \eqref{eq:sample_mean_mgf}
follows from independence and
\eqref{eq:mgf_representation}.
\end{pf}

Since \(M_n(z,x)\) is finite in a neighborhood of \(z=0\),
the distribution \(Q_{n,x}\) has moments of all orders. For
a Borel measurable function \(f\colon\mathbb R\to\mathbb R\) of at most
polynomial growth, define the positive linear operator
\[
(L_nf)(x)
=
\mathbb E_xf(\eta_{n,x})
=
\int_{\mathbb R}f(t)\,Q_{n,x}(dt),
\qquad x\in\mathcal M.
\]
Applying \eqref{eq:weak_differential_identity} to the family
\(\{Q_{n,x}\}\), whose variance function is \(V/n\), gives
\[
\frac{V(x)}{n}(L_nf)'(x)
=
\int_{\mathbb R}(t-x)f(t)\,Q_{n,x}(dt),
\qquad x\in\mathcal M.
\]
Thus, the sample-mean operators satisfy the differential
identity characteristic of operators of exponential type;
see \cite{May1976,IsmailMay1978}.

For \(e_k(t)=t^k\), \(k=0,1,2\), we have
\[
L_ne_0=1,
\qquad
(L_ne_1)(x)=x,
\qquad
(L_ne_2)(x)=x^2+\frac{V(x)}{n},
\]
and, in particular,
\[
\int_{\mathbb R}(t-x)^2\,Q_{n,x}(dt)
=
\frac{V(x)}{n}.
\]

\begin{proposition}
\label{prop:operator_convergence}
For every function \(f\in C_b(\mathbb R)\) and every compact
set \(K\subset\mathcal M\),
\begin{equation}
\label{eq:operator_uniform_convergence}
\lim_{n\to\infty}
\sup_{x\in K}|(L_nf)(x)-f(x)|
=
0.
\end{equation}
\end{proposition}

\begin{pf}
Let
\[
K_1
=
\{t\in\mathbb R:\operatorname{dist}(t,K)\leq1\},
\qquad
V_K=\sup_{x\in K}V(x).
\]
Denote by \(\omega_{K_1}(f,\delta)\) the modulus of
continuity of \(f\) on \(K_1\). For \(0<\delta\leq1\), the
uniform continuity of \(f\) on \(K_1\), the decomposition
according to the events
\(|\eta_{n,x}-x|\leq\delta\) and
\(|\eta_{n,x}-x|>\delta\), and Chebyshev's inequality give
\[
\sup_{x\in K}|(L_nf)(x)-f(x)|
\leq
\omega_{K_1}(f,\delta)
+
\frac{2\|f\|_\infty V_K}{n\delta^2}.
\]
First letting \(n\to\infty\) and then
\(\delta\downarrow0\), we obtain
\eqref{eq:operator_uniform_convergence}.
\end{pf}

Thus, the operators \(L_n\) generated by sample means
approximate the identity operator locally uniformly on
\(\mathcal M\). Formula \eqref{eq:sample_mean_mgf} will be
used below to compute cumulants, central moments, and
asymptotic expansions.

\section{Cumulants, Central Moments, and Coefficient Operators}
\label{sec:cumulants-and-moments}

In this section, the description of the family in terms of
its variance function is translated into the quantities that
determine the asymptotic expansion of the sample-mean
operators. We first derive the cumulants of the original
family, then pass to the cumulants and central moments of the
sample means, and finally collect the coefficients of equal
powers of \(n^{-1}\) into differential operators
\(\mathcal A_r[V]\).

\subsection{Cumulant Recursion}
\label{subsec:cumulant-recursion}

Let
\(
C(z,x)=\log M(z,x)
\)
be the cumulant generating function of the distribution
\(Q_x\). For every real \(z\) for which \(M(z,x)\) is finite,
we have \(M(z,x)>0\), so the function \(C\) is well defined.
The canonical representation of the family gives
\[
C(z,x)
=
K\bigl(\theta(x)+z\bigr)-K\bigl(\theta(x)\bigr).
\]
Thus, \(C\) is real-analytic in a neighborhood of \(z=0\),
where
\(
C(z,x)
=
\sum_{j=1}^{\infty}
\kappa_j(x)\frac{z^j}{j!},
\)
and
\(
\kappa_j(x)
=
\left.
\frac{\partial^j C}{\partial z^j}(z,x)
\right|_{z=0}
\)
is the \(j\)-th cumulant of the distribution \(Q_x\).

\begin{proposition}
\label{prop:cumulant-recursion}
The cumulants of the family satisfy the relations
\[
\kappa_1(x)=x,
\qquad
\kappa_{j+1}(x)
=
V(x)\frac{d}{dx}\kappa_j(x),
\qquad j\geq1.
\]
In particular,
\begin{equation}
\label{eq:cumulant-operator-form}
\kappa_j(x)
=
\left(
V(x)\frac{d}{dx}
\right)^{j-2}V(x),
\qquad j\geq2,
\end{equation}
where the zeroth power denotes the identity operator.
\end{proposition}

\begin{pf}
Dividing the differential identity for \(M\) by \(M>0\),
we obtain
\begin{equation}
\label{eq:cumulant-pde}
\frac{\partial C}{\partial z}(z,x)
-
V(x)\frac{\partial C}{\partial x}(z,x)
=
x.
\end{equation}
Since \(C(0,x)=0\), setting \(z=0\) gives
\(\kappa_1(x)=x\). Differentiating
\eqref{eq:cumulant-pde} \(j\) times with respect to \(z\) and
then setting \(z=0\), we obtain
\[
\kappa_{j+1}(x)
=
V(x)\frac{d}{dx}\kappa_j(x),
\qquad j\geq1.
\]
Successive application of this recursion, starting with
\(\kappa_2=V\), proves
\eqref{eq:cumulant-operator-form}.
\end{pf}

Thus, the first cumulants are
\[
\begin{aligned}
\kappa_1(x)&=x, \quad
\kappa_2(x)=V(x), \quad
\kappa_3(x)=V(x)V'(x),\quad
\kappa_4(x)=
V(x)\left[
\bigl(V'(x)\bigr)^2+V(x)V''(x)
\right].
\end{aligned}
\]
Consequently, each cumulant of order \(j\geq2\) is a
differential polynomial in \(V\) and its derivatives up to
order \(j-2\).

For the following subclass of quadratic variance functions,
this recursion can be solved explicitly.

\begin{corollary}
\label{cor:quadratic-variance-cumulants}
Suppose that
\(
V(x)=x(\alpha x+\beta),
\)
\(
x\in\mathcal M,
\) 
for some pair \((\alpha,\beta)\in\mathbb R^2\). Then, for
every \(k\geq1\),
\[
\kappa_{k+1}(x)
=
x(\alpha x+\beta)
\sum_{m=1}^{k}
m!\,\beta^{\,k-m}S(k,m)(\alpha x)^{m-1},
\]
where \(S(k,m)\) denotes a Stirling number of the second kind.
Factors with exponent zero are understood to be equal to
one.
\end{corollary}

\begin{pf}
Define
\[
P_k(u)
=
\sum_{m=1}^{k}
m!\,\beta^{\,k-m}S(k,m)u^{m-1}.
\]
The recurrence relation for the Stirling numbers
\[
S(k+1,m)=mS(k,m)+S(k,m-1)
\]
gives
\[
P_{k+1}(u)
=
(\beta+2u)P_k(u)
+
u(u+\beta)P_k'(u).
\]
Since \(P_1=1\), the required formula holds for
\(\kappa_2=V\). If
\(\kappa_{k+1}(x)=V(x)P_k(\alpha x)\), then, with
\(u=\alpha x\), the cumulant recursion gives
\[
\begin{aligned}
\kappa_{k+2}(x)
&=
V(x)\frac{d}{dx}
\left[V(x)P_k(\alpha x)\right]
=
V(x)
\left[
(\beta+2u)P_k(u)
+
u(u+\beta)P_k'(u)
\right]
=
V(x)P_{k+1}(u).
\end{aligned}
\]
The assertion follows by induction.
\end{pf}

The Bernoulli, Poisson, geometric, and exponential natural
exponential families correspond, respectively, to
\[
(\alpha,\beta)
=
(-1,1),\qquad
(0,1),\qquad
(1,1),\qquad
(1,0).
\]
These examples belong to the class of natural exponential
families with quadratic variance functions classified by
Morris~\cite{Morris1982}.

\subsection{Cumulants of Sample Means}
\label{subsec:sample-mean-cumulants}

Let
\(
C_n(z,x)=\log M_n(z,x)
\)
be the cumulant generating function of the sample mean
\(\eta_{n,x}\), and denote
\[
\kappa_j^{(n)}(x)
=
\left.
\frac{\partial^j C_n}{\partial z^j}(z,x)
\right|_{z=0}.
\]
The sample-mean identity
\eqref{eq:sample_mean_mgf} gives
\begin{equation}
\label{eq:sample-mean-cgf}
C_n(z,x)
=
nC\left(\frac{z}{n},x\right)
\end{equation}
in a neighborhood of \(z=0\).

\begin{proposition}
\label{prop:sample-mean-cumulants}
For all \(n\geq1\), \(x\in\mathcal M\), and \(j\geq1\),
\[
\kappa_j^{(n)}(x)
=
\frac{\kappa_j(x)}{n^{j-1}}.
\]
\end{proposition}

\begin{pf}
Differentiate \eqref{eq:sample-mean-cgf} \(j\) times with
respect to \(z\) and then set \(z=0\).
\end{pf}

Set
\(
\Delta_{n,x}=\eta_{n,x}-x.
\)
Since \(\kappa_1(x)=x\), the centered cumulant generating
function takes the form
\begin{equation}
\label{eq:centered-sample-mean-cgf}
\log\mathbb E_x e^{z\Delta_{n,x}}
=
\sum_{j=2}^{\infty}
\frac{\kappa_j(x)}{n^{j-1}}
\frac{z^j}{j!}.
\end{equation}
Thus, centering removes the first cumulant, while the
\(j\)-th of the remaining cumulants is of order
\(n^{-(j-1)}\). This scaling determines the powers of
\(n^{-1}\) that occur in the central moments and,
consequently, in the asymptotic coefficient operators.

\subsection{Central Moments and Bell Polynomials}
\label{subsec:central-moments}

For \(k\geq0\), set
\(
\mu_k^{(n)}(x)
=
\mathbb E_x\Delta_{n,x}^k.
\)
Then
\(
\mu_0^{(n)}(x)=1,
\)
\(
\mu_1^{(n)}(x)=0.
\)

For \(1\leq q\leq k\), denote by
\[
B_{k,q}(y_1,\ldots,y_{k-q+1})
=
k!
\sum_{\substack{
m_1,\ldots,m_{k-q+1}\geq0\\
\sum_{j=1}^{k-q+1}m_j=q\\
\sum_{j=1}^{k-q+1}jm_j=k
}}
\prod_{j=1}^{k-q+1}
\frac{1}{m_j!}
\left(\frac{y_j}{j!}\right)^{m_j}.
\]
the partial exponential Bell polynomial. For this notation
and the standard moment--cumulant formula used below, see,
for example, \cite{Comtet1974}.

\begin{proposition}
\label{prop:central-moments}
For all \(n\geq1\), \(x\in\mathcal M\), and \(k\geq2\),
\begin{equation}
\label{eq:central-moment-expansion}
\mu_k^{(n)}(x)
=
\sum_{r=\lceil k/2\rceil}^{k-1}
\frac{1}{n^r}
B_{k,k-r}
\bigl(
0,\kappa_2(x),\ldots,\kappa_{r+1}(x)
\bigr).
\end{equation}
\end{proposition}

\begin{pf}
The moment--cumulant formula together with
\eqref{eq:centered-sample-mean-cgf} gives
\[
\mu_k^{(n)}(x)
=
\sum_{q=1}^{\lfloor k/2\rfloor}
B_{k,q}
\left(
0,
\frac{\kappa_2(x)}{n},
\ldots,
\frac{\kappa_{k-q+1}(x)}{n^{k-q}}
\right).
\]
Since the first argument of the Bell polynomial is zero, we
have \(m_1=0\) for every nonzero monomial. Therefore, the
power of \(n^{-1}\) in such a monomial is
\[
\sum_{j=2}^{k-q+1}(j-1)m_j
=
k-q.
\]
Consequently,
\[
\mu_k^{(n)}(x)
=
\sum_{q=1}^{\lfloor k/2\rfloor}
\frac{1}{n^{k-q}}
B_{k,q}
\bigl(
0,\kappa_2(x),\ldots,\kappa_{k-q+1}(x)
\bigr).
\]
The substitution \(r=k-q\) yields
\eqref{eq:central-moment-expansion}.
\end{pf}

In particular,
\[
\begin{aligned}
\mu_2^{(n)}(x)
&= \frac{V(x)}{n}, \quad
\mu_3^{(n)}(x) = \frac{V(x)V'(x)}{n^2}, \\
\mu_4^{(n)}(x)
&= \frac{3V(x)^2}{n^2}
+ \frac{V(x)\left[\bigl(V'(x)\bigr)^2+V(x)V''(x)\right]}{n^3}.
\end{aligned}
\]

For a fixed \(r\geq1\), a term of order \(n^{-r}\) can occur
only in moments with indices satisfying
\begin{equation}
\label{eq:moment-index-range}
r+1\leq k\leq2r.
\end{equation}
This is precisely the range of derivative orders that enter
the coefficient operator \(\mathcal A_r[V]\).

\begin{corollary}
\label{cor:uniform-even-moments}
For every integer \(q\geq1\) and every compact set
\(J\subset\mathcal M\),
\[
\sup_{x\in J}
\mathbb E_x
\left|\eta_{n,x}-x\right|^{2q}
=
O(n^{-q}),
\qquad n\to\infty.
\]
\end{corollary}

\begin{pf}
For \(k=2q\), the expansion
\eqref{eq:central-moment-expansion} begins with the power
\(n^{-q}\). Its finitely many coefficient functions are
continuous and hence bounded on \(J\). Since an even central
moment coincides with the corresponding absolute moment, the
assertion follows.
\end{pf}

\subsection{Coefficient Operators}
\label{subsec:coefficient-operators}

We now collect the coefficients of equal powers of \(n^{-1}\)
that arise when the Taylor expansion of \(f\) is averaged
with respect to the distribution \(Q_{n,x}\). According to
\eqref{eq:central-moment-expansion} and
\eqref{eq:moment-index-range}, the coefficient of \(n^{-r}\)
can involve only the derivatives
\(
f^{(k)},
\) 
\(
r+1\leq k\leq2r.
\)

\begin{definition}
\label{def:coefficient-operators}
Set
\(
\mathcal A_0[V]f=f.
\)
For \(r\geq1\) and \(f\in C^{2r}(\mathcal M)\), define
\begin{equation}
\label{eq:coefficient-operator}
\mathcal A_r[V]f(x)
=
\sum_{k=r+1}^{2r}
\frac{
B_{k,k-r}
\bigl(
0,\kappa_2(x),\ldots,\kappa_{r+1}(x)
\bigr)
}{k!}
f^{(k)}(x).
\end{equation}
\end{definition}

Substitution of the multi-index definition of the Bell
polynomials gives the equivalent representation
\begin{equation}
\label{eq:coefficient-operator-multi-index}
\mathcal A_r[V]f(x)
=
\sum_{\substack{
m_2,\ldots,m_{r+1}\geq0\\
\sum_{j=2}^{r+1}(j-1)m_j=r
}}
\left(
\prod_{j=2}^{r+1}
\frac{\kappa_j(x)^{m_j}}
     {m_j!\,(j!)^{m_j}}
\right)
f^{\left(\sum_{j=2}^{r+1}jm_j\right)}(x).
\end{equation}
Indeed, for every admissible multi-index,
\(
k=\sum_{j=2}^{r+1}jm_j,
\) 
\(
k-r=\sum_{j=2}^{r+1}m_j,
\)
which coincides exactly with the relation between the
indices in \(B_{k,k-r}\).

The same representation determines the principal part of
the operator. For every admissible multi-index,
\[
\sum_{j=2}^{r+1}jm_j
=
r+\sum_{j=2}^{r+1}m_j
\leq2r.
\]
Equality can occur only when
\[
m_2=r,
\qquad
m_j=0,
\quad 
3\leq j\leq r+1.
\]
Since \(\kappa_2=V\), this unique multi-index contributes the
term
\(
\frac{V^r}{r!\,(2!)^r}f^{(2r)}.
\)
Therefore,
\begin{equation}
\label{eq:coefficient-operator-principal-part}
\mathcal A_r[V]f(x)
=
\frac{V(x)^r}{2^r r!}f^{(2r)}(x)
+
\sum_{k=r+1}^{2r-1}
a_{r,k}(x)f^{(k)}(x),
\qquad r\geq1,
\end{equation}
where \(a_{r,k}\in C^\infty(\mathcal M)\), and an empty sum
is understood to be zero.

According to \eqref{eq:cumulant-operator-form}, all
coefficients of the operator \(\mathcal A_r[V]\) are
differential polynomials in \(V\). They involve derivatives
of \(V\) of order at most \(r-1\), whereas the operator
itself has exact differential order \(2r\), since \(V>0\) on
\(\mathcal M\).

The first coefficient operators are
\begin{equation}
\label{eq:first-coefficient-operators}
\begin{aligned}
\mathcal A_1[V]f(x)
&=
\frac{V(x)}{2}f''(x),\\[1mm]
\mathcal A_2[V]f(x)
&=
\frac{V(x)V'(x)}{6}f'''(x)
+
\frac{V(x)^2}{8}f^{(4)}(x),\\[1mm]
\mathcal A_3[V]f(x)
&=
\frac{
V(x)\left[
\bigl(V'(x)\bigr)^2+V(x)V''(x)
\right]
}{24}
f^{(4)}(x)
+
\frac{V(x)^2V'(x)}{12}f^{(5)}(x)
+
\frac{V(x)^3}{48}f^{(6)}(x).
\end{aligned}
\end{equation}

Thus, the variance function determines the entire hierarchy
of cumulants, central moments, and coefficient operators
that govern the asymptotic expansion of \(L_n\) to arbitrary
order.

\section{Arbitrary-Order Expansions for the Operators
  \texorpdfstring{\(L_n\)}{L\_n}}
\label{sec:all-order-expansion}
The coefficient operators introduced in the preceding
section describe the ordinary, that is, undifferentiated,
expansion of \(L_nf\). We first establish a quantitative
estimate for an arbitrary integer order of smoothness and
then derive an asymptotic expansion of arbitrary order and
the corresponding Voronovskaya-type formulas.

General asymptotic expansions of positive linear operators
in terms of central moments were developed by
Sikkema~\cite{Sikkema1970}. Among later complete expansions
for individual classical operators, we mention the result of
Tachev for Bernstein operators~\cite{Tachev2012}.
Complete pointwise expansions were also obtained for several
specific operators of exponential type associated with
nonquadratic functions \(p\); see
\cite{AbelGupta2021,AbelGuptaKushnirevych2021,
AbelKushnirevych2025}.
In our setting, the exact moment identities from the
preceding section lead to the following quantitative result,
uniform on compact subsets, with all its coefficients
determined by the variance function \(V\).

\subsection{Quantitative Expansion of Arbitrary Order}
\label{subsec:main-expansion}

\begin{theorem}[Quantitative Expansion of Arbitrary Order]
\label{thm:all-order-expansion}
Let \(m\geq1\) be an integer, and set
\(
s=\left\lfloor\frac{m}{2}\right\rfloor.
\)
Let \(f\colon\mathbb R\to\mathbb R\) be a Borel measurable function of
at most polynomial growth, that is, for some \(c>0\) and
\(p\geq0\),
\[
|f(t)|
\leq
c\bigl(1+|t|^p\bigr),
\qquad t\in\mathbb R.
\]
Suppose that
\(
f|_{\mathcal M}\in C^m(\mathcal M).
\)

Let \(J\subset\mathcal M\) be compact, and choose \(\rho>0\)
so that
\[
J_\rho
=
\bigl\{
y\in\mathbb R:
\operatorname{dist}(y,J)\leq\rho
\bigr\}
\subset\mathcal M.
\]
Define
\[
\omega_{m,J_\rho}(f;\delta)
=
\sup_{\substack{
u,v\in J_\rho\\
|u-v|\leq\delta
}}
\left|
f^{(m)}(u)-f^{(m)}(v)
\right|.
\]
Then there exist constants \(C,D>0\), independent of \(n\),
such that, for all sufficiently large \(n\),
\begin{equation}
\label{eq:quantitative-remainder}
\begin{aligned}
\sup_{x\in J}
\left|
L_nf(x)
-
\sum_{r=0}^{s}
\frac{1}{n^r}\mathcal A_r[V]f(x)
\right|
\leq{}&
C n^{-m/2}
\omega_{m,J_\rho}\bigl(f;n^{-1/2}\bigr)
+
D n^{-s-1}.
\end{aligned}
\end{equation}
\end{theorem}

\begin{pf}
Set
\(
\Delta_{n,x}=\eta_{n,x}-x
\)
and
\[
R_m(x,h)
=
f(x+h)
-
\sum_{k=0}^{m}
\frac{f^{(k)}(x)}{k!}h^k.
\]
If \(x\in J\) and \(|h|\leq\rho\), then the line segment
joining \(x\) and \(x+h\) is contained in \(J_\rho\).
Taylor's formula, followed by partitioning this segment into
subintervals of length at most \(\delta\), gives, for
\(0<\delta\leq\rho\),
\begin{equation}
\label{eq:taylor-modulus-remainder}
|R_m(x,h)|
\leq
\frac{|h|^m}{m!}
\left(
1+\frac{|h|}{\delta}
\right)
\omega_{m,J_\rho}(f;\delta).
\end{equation}

Corollary~\ref{cor:uniform-even-moments} and Lyapunov's
inequality imply that, for every integer \(k\geq0\),
\[
\sup_{x\in J}
\mathbb E_x|\Delta_{n,x}|^k
=
O\left(n^{-k/2}\right).
\]
Setting \(\delta=n^{-1/2}\) in
\eqref{eq:taylor-modulus-remainder}, we obtain
\[
\begin{aligned}
\sup_{x\in J}
\mathbb E_x\bigl[
&|R_m(x,\Delta_{n,x})|;
|\Delta_{n,x}|\leq\rho
\bigr]
\leq
C n^{-m/2}
\omega_{m,J_\rho}\bigl(f;n^{-1/2}\bigr).
\end{aligned}
\]

On the complementary event, the polynomial growth
assumption and the boundedness on \(J\) of the derivatives
\(f^{(k)}\), \(0\leq k\leq m\), imply that
\[
|R_m(x,h)|
\leq
C_J\bigl(1+|h|^d\bigr),
\qquad
x\in J,\qquad h\in\mathbb R,
\]
where \(d=\max\{p,m\}\). Choose an integer \(q\geq s+1\)
such that \(2q\geq d\). For \(|h|>\rho\),
\[
1+|h|^d\leq C_{\rho,q}|h|^{2q}.
\]
Another application of
Corollary~\ref{cor:uniform-even-moments} gives
\[
\sup_{x\in J}
\mathbb E_x\bigl[
|R_m(x,\Delta_{n,x})|;
|\Delta_{n,x}|>\rho
\bigr]
=
O\left(n^{-s-1}\right).
\]

Combining the local and tail estimates, we obtain, uniformly
for \(x\in J\),
\begin{equation}
\label{eq:quantitative-taylor-expansion}
\begin{aligned}
L_nf(x)
={}&
\sum_{k=0}^{m}
\frac{f^{(k)}(x)}{k!}\mu_k^{(n)}(x)
+
O\left(
n^{-m/2}
\omega_{m,J_\rho}\bigl(f;n^{-1/2}\bigr)
\right)
+
O\left(n^{-s-1}\right).
\end{aligned}
\end{equation}
By \eqref{eq:central-moment-expansion} and
\eqref{eq:moment-index-range}, the terms of orders
\(n^{-r}\), \(0\leq r\leq s\), in the finite sum in
\eqref{eq:quantitative-taylor-expansion} form precisely
\(\mathcal A_r[V]f\). All remaining terms are of order
\(O(n^{-s-1})\). This proves
\eqref{eq:quantitative-remainder}.
\end{pf}

\begin{corollary}[Asymptotic Expansion of Arbitrary Order]
\label{cor:ordinary-all-order-expansion}
Let \(N\geq1\), and let
\(f\colon\mathbb R\to\mathbb R\) be a Borel measurable function of at
most polynomial growth such that
\(
f|_{\mathcal M}\in C^{2N}(\mathcal M).
\)
Then, for every compact set \(J\subset\mathcal M\),
\begin{equation}
\label{eq:all-order-expansion}
\sup_{x\in J}
\left|
L_nf(x)
-
\sum_{r=0}^{N}
\frac{1}{n^r}\mathcal A_r[V]f(x)
\right|
=
o\left(n^{-N}\right).
\end{equation}

If, in addition,
\[
f|_{\mathcal M}\in C^{2N+2}(\mathcal M),
\]
then
\begin{equation}
\label{eq:refined-all-order-expansion}
\sup_{x\in J}
\left|
L_nf(x)
-
\sum_{r=0}^{N+1}
\frac{1}{n^r}\mathcal A_r[V]f(x)
\right|
=
o\left(n^{-N-1}\right),
\end{equation}
and consequently,
\begin{equation}
\label{eq:all-order-expansion-big-O}
\sup_{x\in J}
\left|
L_nf(x)
-
\sum_{r=0}^{N}
\frac{1}{n^r}\mathcal A_r[V]f(x)
\right|
=
O\left(n^{-N-1}\right).
\end{equation}
\end{corollary}

\begin{pf}
Apply Theorem~\ref{thm:all-order-expansion} with \(m=2N\).
Since \(f^{(2N)}\) is uniformly continuous on every compact
neighborhood \(J_\rho\Subset\mathcal M\),
\[
\omega_{2N,J_\rho}(f;\delta)\longrightarrow0,
\qquad \delta\downarrow0,
\]
which gives \eqref{eq:all-order-expansion}. Applying the same
theorem with \(m=2N+2\) gives
\eqref{eq:refined-all-order-expansion}. Finally,
\(\mathcal A_{N+1}[V]f\) is bounded on \(J\), so separating
the term with \(r=N+1\) gives
\eqref{eq:all-order-expansion-big-O}.
\end{pf}

\begin{rmk}[Intermediate Smoothness]
\label{rem:intermediate-ordinary-rates}
If \(m=2N+1\), then
Theorem~\ref{thm:all-order-expansion} gives
\[
\sup_{x\in J}
\left|
L_nf(x)
-
\sum_{r=0}^{N}
\frac{1}{n^r}\mathcal A_r[V]f(x)
\right|
=
o\left(n^{-N-1/2}\right).
\]
More generally, if, for some \(0<\gamma\leq1\),
\[
\omega_{m,J_\rho}(f;\delta)
=
O\left(\delta^\gamma\right),
\qquad \delta\downarrow0,
\]
then
\[
\sup_{x\in J}
\left|
L_nf(x)
-
\sum_{r=0}^{s}
\frac{1}{n^r}\mathcal A_r[V]f(x)
\right|
=
O\left(n^{-(m+\gamma)/2}\right).
\]
\end{rmk}

\subsection{Higher-Order Voronovskaya Formulas}
\label{subsec:voronovskaya-formulas}

The expansion of arbitrary order immediately yields
Voronovskaya-type formulas of arbitrary order.

\begin{corollary}
\label{cor:higher-order-voronovskaya}
Let $N\geq1$, and let
$f\colon\mathbb R\to\mathbb R$ be a Borel measurable function
of at most polynomial growth such that
\(
f|_{\mathcal M}\in C^{2N}(\mathcal M).
\)
Then, for every compact set $J\subset\mathcal M$,
\[
\lim_{n\to\infty}
\sup_{x\in J}
\left|
n^N
\left[
L_nf(x)
-
\sum_{r=0}^{N-1}
\frac{1}{n^r}\mathcal A_r[V]f(x)
\right]
-
\mathcal A_N[V]f(x)
\right|
=
0.
\]
\end{corollary}

\begin{pf}
Separate the term corresponding to $r=N$ in
\eqref{eq:all-order-expansion} and multiply the resulting
relation by $n^N$.
\end{pf}

For $N=1$, if $f|_{\mathcal M}\in C^2(\mathcal M)$, we
obtain the classical first-order Voronovskaya
formula~\cite{Voronovskaya1932}:
\[
n\bigl(L_nf(x)-f(x)\bigr)
\longrightarrow
\frac{V(x)}{2}f''(x)
\]
locally uniformly on $\mathcal M$.

For $N=2$, if $f|_{\mathcal M}\in C^4(\mathcal M)$, then
\[
\begin{aligned}
n^2
\left[
L_nf(x)
-
f(x)
-
\frac{V(x)}{2n}f''(x)
\right]
\longrightarrow{}&
\frac{V(x)V'(x)}{6}f'''(x)
+
\frac{V(x)^2}{8}f^{(4)}(x)
\end{aligned}
\]
locally uniformly on $\mathcal M$.

Thus, the operators $L_n$ admit a complete hierarchy of
Voronovskaya-type formulas whose coefficients are determined
entirely by the variance function $V$.

\section{Quantitative Simultaneous Expansions}
\label{sec:simultaneous}

In this section, we establish quantitative asymptotic
expansions for the derivatives of the operators $L_n$,
uniformly on compact subsets. General Voronovskaya-type and
asymptotic theorems in simultaneous approximation were
obtained by Gonska and P\u{a}lt\u{a}nea
\cite{GonskaPaltanea2010}; related results for derivatives of
positive linear operators were developed by
Holho\c{s}~\cite{Holhos2019,Holhos2022Simultaneous}. For
operators of exponential type, see
Holho\c{s}~\cite{Holhos2022Exponential}, while
Abel~\cite{Abel2023} established that full asymptotic
expansions can be differentiated term by term pointwise. In
the present framework of sample-mean operators, the
structure of the exponential family leads to the
quantitative estimate below, uniform on compact subsets,
with coefficient operators explicitly determined by the
variance function $V$.

Recall that
\[
(L_nf)(x)
=
\int_{\mathbb R}f(t)\,Q_{n,x}(dt),
\qquad x\in\mathcal M,
\]
where
\[
Q_{n,x}(dt)
=
\exp\{n\tau(x)t-nH(x)\}\,\nu_n(dt),
\]
the measure $\nu_n$ is independent of $x$, and
\[
\tau'(x)=\frac{1}{V(x)},
\qquad
H'(x)=\frac{x}{V(x)}.
\]
Hence,
\[
\frac{\partial}{\partial x}
\bigl(n\tau(x)t-nH(x)\bigr)
=
\frac{n(t-x)}{V(x)}.
\]

For $n\geq1$, recursively define the logarithmic-derivative
polynomials of the distribution (the score polynomials) by
\begin{equation}
\label{eq:score-polynomial-recursion}
P_{0,n}(t,x)=1,
\qquad
P_{\ell+1,n}(t,x)
=
\frac{\partial}{\partial x}P_{\ell,n}(t,x)
+
\frac{n(t-x)}{V(x)}P_{\ell,n}(t,x),
\end{equation}
where the partial derivative is taken with $t$ fixed.

\subsection{Score Polynomials and Weighted Moments}

\begin{lemma}
\label{lem:score-representation}
Let $f\colon\mathbb R\to\mathbb R$ be a Borel measurable
function of at most polynomial growth. Then, for every
$n\geq1$,
\[
L_nf\in C^\infty(\mathcal M),
\]
and
\begin{equation}
\label{eq:score-representation}
(L_nf)^{(\ell)}(x)
=
\int_{\mathbb R}
f(t)P_{\ell,n}(t,x)\,Q_{n,x}(dt),
\qquad
\ell\in\mathbb N_0,\qquad x\in\mathcal M.
\end{equation}
\end{lemma}

\begin{pf}
Fix $n\geq1$ and a compact set $K\subset\mathcal M$.
Choose points in $\mathcal M$ such that
\[
a_-<a<\min K\leq\max K<b<b_+.
\]
Since $\tau$ is strictly increasing, for $x\in[a,b]$ we have
\[
e^{n\tau(x)t}
\leq
\begin{cases}
e^{n\tau(b)t},&t\geq0,\\
e^{n\tau(a)t},&t<0.
\end{cases}
\]

Let $R(t,x)$ be a polynomial in $t$ with coefficients in
$C^1(\mathcal M)$. On $[a,b]$, the coefficient functions of
$R$ and $\partial_xR$, as well as $V^{-1}$ and $e^{-nH}$,
are bounded. Since $f$ is of at most polynomial growth, the
integrands arising from
\[
f(t)R(t,x),
\qquad
f(t)\frac{\partial R}{\partial x}(t,x),
\qquad
f(t)(t-x)R(t,x)
\]
are, uniformly in $x\in[a,b]$, dominated by a constant
multiple of
\[
(1+|t|^q)
\left(
e^{n\tau(b)t}\mathbf 1_{\{t\geq0\}}
+
e^{n\tau(a)t}\mathbf 1_{\{t<0\}}
\right)
\]
for some $q\geq0$. The polynomial factor can be absorbed
into an exponential with a slightly larger exponent:
\[
\begin{aligned}
(1+|t|^q)e^{n\tau(b)t}\mathbf 1_{\{t\geq0\}}
&\leq
C e^{n\tau(b_+)t}\mathbf 1_{\{t\geq0\}},\\
(1+|t|^q)e^{n\tau(a)t}\mathbf 1_{\{t<0\}}
&\leq
C e^{n\tau(a_-)t}\mathbf 1_{\{t<0\}}.
\end{aligned}
\]
Both majorants are integrable with respect to $\nu_n$,
since the normalization of $Q_{n,y}$ implies that
\[
\int_{\mathbb R}e^{n\tau(y)t}\,\nu_n(dt)
=
e^{nH(y)}<\infty,
\qquad y\in\mathcal M.
\]

Therefore, differentiation under the integral sign gives
\[
\begin{aligned}
\frac{d}{dx}
\int_{\mathbb R}f(t)R(t,x)\,Q_{n,x}(dt)
=
\int_{\mathbb R}
f(t)
\left[
\frac{\partial R}{\partial x}(t,x)
+
\frac{n(t-x)}{V(x)}R(t,x)
\right]
Q_{n,x}(dt).
\end{aligned}
\]
Starting with $R=P_{0,n}=1$ and using
\eqref{eq:score-polynomial-recursion}, we obtain
\eqref{eq:score-representation} by induction. Since the
compact set $K$ was arbitrary, it follows that
$L_nf\in C^\infty(\mathcal M)$.
\end{pf}

\begin{lemma}
\label{lem:score-structure}
For every $\ell\in\mathbb N_0$, there exist a finite set
$\mathcal I_\ell\subset\mathbb N_0^2$ and functions
\(
c_{\ell,a,b}\in C^\infty(\mathcal M),
\) 
\((a,b)\in\mathcal I_\ell,
\)
independent of $n$, such that
\[
P_{\ell,n}(t,x)
=
\sum_{(a,b)\in\mathcal I_\ell}
c_{\ell,a,b}(x)n^a(t-x)^b.
\]
Every pair $(a,b)\in\mathcal I_\ell$ satisfies
\[
0\leq a\leq\ell,
\qquad
0\leq b\leq\ell,
\qquad
2a-b\leq\ell.
\]
\end{lemma}

\begin{pf}
For $\ell=0$, the assertion follows from $P_{0,n}=1$.
Suppose that it holds for some $\ell\geq0$, and consider a
term
\(
c(x)n^a(t-x)^b.
\)
Differentiation with respect to $x$ and multiplication by
the score factor in \eqref{eq:score-polynomial-recursion}
produce terms with indices
\[
(a,b),
\qquad
(a,b-1),
\qquad
(a+1,b+1),
\]
respectively, where the term with index $(a,b-1)$ is absent
when $b=0$. If $w=2a-b$, then the corresponding weights are
$w$, $w+1$, and $w+1$. Thus, $w\leq\ell$ implies that each
new weight does not exceed $\ell+1$, while the bounds on $a$
and $b$ remain valid after replacing $\ell$ by $\ell+1$.
The new coefficients are obtained by differentiation and
multiplication by $V^{-1}$; they belong to
$C^\infty(\mathcal M)$ because $V$ is positive and
real-analytic on $\mathcal M$.
\end{pf}

\begin{lemma}
\label{lem:weighted-moment}
Let $J\subset\mathcal M$ be compact. For every
$\ell\in\mathbb N_0$ and every real number $k\geq0$, there
exists a constant $C_{\ell,k,J}>0$ such that
\begin{equation}
\label{eq:weighted-moment}
\sup_{x\in J}
\mathbb E_x
\left[
|P_{\ell,n}(\eta_{n,x},x)|
|\eta_{n,x}-x|^k
\right]
\leq
C_{\ell,k,J}n^{(\ell-k)/2},
\qquad n\geq1.
\end{equation}
\end{lemma}

\begin{pf}
Set $\Delta_{n,x}=\eta_{n,x}-x$. The previously established
estimate for even moments, uniform on compact subsets, gives,
for every integer $q\geq1$,
\[
\sup_{x\in J}\mathbb E_x|\Delta_{n,x}|^{2q}
\leq C_{q,J}n^{-q}.
\]
Therefore, Lyapunov's inequality implies that, for every real
number $u\geq0$,
\[
\sup_{x\in J}\mathbb E_x|\Delta_{n,x}|^u
\leq C_{u,J}n^{-u/2}.
\]
By Lemma~\ref{lem:score-structure}, the coefficient functions
$c_{\ell,a,b}$ are bounded on $J$. Hence,
\[
\begin{aligned}
\mathbb E_x
\left[
|P_{\ell,n}(\eta_{n,x},x)|
|\Delta_{n,x}|^k
\right]
&\leq
C\sum_{(a,b)\in\mathcal I_\ell}
n^a\mathbb E_x|\Delta_{n,x}|^{b+k}
\leq
C\sum_{(a,b)\in\mathcal I_\ell}
n^{a-(b+k)/2}.
\end{aligned}
\]
Since $2a-b\leq\ell$,
\[
a-\frac{b+k}{2}
\leq
\frac{\ell-k}{2}.
\]
The set $\mathcal I_\ell$ is finite, and therefore
\eqref{eq:weighted-moment} follows.
\end{pf}

\subsection{The Main Simultaneous Expansion}

For a compact set $J\subset\mathcal M$ and $\rho>0$ such
that
\(
J_\rho
=
\{y\in\mathbb R:\operatorname{dist}(y,J)\leq\rho\}
\Subset\mathcal M,
\)
define
\[
\omega_{m,J_\rho}(f;\delta)
=
\sup_{\substack{u,v\in J_\rho\\|u-v|\leq\delta}}
|f^{(m)}(u)-f^{(m)}(v)|.
\]

\begin{theorem}[Quantitative Simultaneous Expansion]
\label{thm:simultaneous-expansion}
Let $\ell,m\in\mathbb N_0$, $m\geq\ell$, and set
\(
s_\ell
=
\left\lfloor\frac{m-\ell}{2}\right\rfloor.
\)
Let $f\colon\mathbb R\to\mathbb R$ be a Borel measurable
function of at most polynomial growth such that
\[
f|_{\mathcal M}\in C^m(\mathcal M).
\]
Then, for every compact set $J\subset\mathcal M$ and every
$\rho>0$ such that $J_\rho\Subset\mathcal M$, there exist
constants $C,D>0$, independent of $n$, such that, for all
sufficiently large $n$,
\begin{equation}
\label{eq:quantitative-simultaneous}
\begin{aligned}
\sup_{x\in J}
\left|
(L_nf)^{(\ell)}(x)
-
\sum_{r=0}^{s_\ell}
\frac{1}{n^r}
\bigl(\mathcal A_r[V]f\bigr)^{(\ell)}(x)
\right|
\leq{}&
C n^{-(m-\ell)/2}
\omega_{m,J_\rho}\bigl(f;n^{-1/2}\bigr)
+
D n^{-s_\ell-1}.
\end{aligned}
\end{equation}
\end{theorem}

For $0\leq r\leq s_\ell$, we have $2r+\ell\leq m$;
therefore, all derivatives occurring in
$\bigl(\mathcal A_r[V]f\bigr)^{(\ell)}$ are well defined
under the assumptions of the theorem.

\begin{pf}
Set
\(
\Delta_{n,x}=\eta_{n,x}-x
\)
and define
\[
R_m(x,h)
=
f(x+h)
-
\sum_{j=0}^{m}\frac{f^{(j)}(x)}{j!}h^j.
\]
By Lemma~\ref{lem:score-representation},
\begin{equation}
\label{eq:simultaneous-basic-decomposition}
\begin{aligned}
(L_nf)^{(\ell)}(x)
={}&
\sum_{j=0}^{m}
\frac{f^{(j)}(x)}{j!}
\mathbb E_x
\left[
P_{\ell,n}(\eta_{n,x},x)\Delta_{n,x}^{\,j}
\right]
+
\mathbb E_x
\left[
P_{\ell,n}(\eta_{n,x},x)
R_m(x,\Delta_{n,x})
\right].
\end{aligned}
\end{equation}

We first estimate the remainder. If $|h|\leq\rho$, then the
segment joining $x$ and $x+h$ is contained in $J_\rho$. For
$m\geq1$, Taylor's formula with integral remainder gives
\[
|R_m(x,h)|
\leq
\frac{|h|^m}{m!}
\sup_{0\leq u\leq1}
|f^{(m)}(x+uh)-f^{(m)}(x)|.
\]
Partitioning the segment into subsegments of length at most
$\delta$, where $0<\delta\leq\rho$, gives
\[
|R_m(x,h)|
\leq
\frac{|h|^m}{m!}
\left(1+\frac{|h|}{\delta}\right)
\omega_{m,J_\rho}(f;\delta).
\]
For $m=0$, the same inequality follows directly by applying
such a partition to $f(x+h)-f(x)$.

Choose $\delta=n^{-1/2}$; for all sufficiently large $n$, we
have $\delta\leq\rho$. Lemma~\ref{lem:weighted-moment},
applied with $k=m$ and $k=m+1$, gives
\[
\begin{aligned}
\sup_{x\in J}\mathbb E_x
\Big[&
|P_{\ell,n}(\eta_{n,x},x)|
|R_m(x,\Delta_{n,x})|;
|\Delta_{n,x}|\leq\rho
\Big]
\leq
C n^{-(m-\ell)/2}
\omega_{m,J_\rho}\bigl(f;n^{-1/2}\bigr).
\end{aligned}
\]

We now consider the event $|\Delta_{n,x}|>\rho$. Let
\(
|f(t)|\leq c(1+|t|^p),
\) 
\(t\in\mathbb R,
\)
and set $d=\max\{p,m\}$. Since the derivatives $f^{(j)}$,
$0\leq j\leq m$, are bounded on $J$, we have
\[
|R_m(x,h)|
\leq
C_J(1+|h|^d),
\qquad x\in J,\qquad h\in\mathbb R.
\]
Choose $u\geq0$ sufficiently large that
\[
\frac{\ell-d-u}{2}\leq-s_\ell-1.
\]
On the event $|\Delta_{n,x}|>\rho$,
\[
1+|\Delta_{n,x}|^d
\leq
C_{\rho,u}|\Delta_{n,x}|^{d+u}.
\]
Another application of Lemma~\ref{lem:weighted-moment} gives
\[
\begin{aligned}
\sup_{x\in J}
\mathbb E_x
\Big[&
|P_{\ell,n}(\eta_{n,x},x)|
|R_m(x,\Delta_{n,x})|;
|\Delta_{n,x}|>\rho
\Big]
\leq
C n^{(\ell-d-u)/2}
\leq
D n^{-s_\ell-1}.
\end{aligned}
\]
Combining the local and tail estimates, we obtain
\begin{equation}
\label{eq:full-score-remainder}
\begin{aligned}
\sup_{x\in J}
\left|
\mathbb E_x
\left[
P_{\ell,n}(\eta_{n,x},x)
R_m(x,\Delta_{n,x})
\right]
\right|
\leq{}&
C n^{-(m-\ell)/2}
\omega_{m,J_\rho}\bigl(f;n^{-1/2}\bigr)
+
D n^{-s_\ell-1}.
\end{aligned}
\end{equation}

It remains to identify the finite sum in
\eqref{eq:simultaneous-basic-decomposition}. If $m=0$, then
$\ell=0$, and this sum is equal to
$f(x)=\mathcal A_0[V]f(x)$. Henceforth, assume that $m\geq1$.

Fix $x_0\in J$ and introduce the Taylor polynomial
\[
p_{x_0}(t)
=
\sum_{j=0}^{m}
\frac{f^{(j)}(x_0)}{j!}(t-x_0)^j.
\]
For every polynomial $p$ of degree at most $m$, Taylor's
formula is exact, and therefore
\[
(L_np)(y)
=
\sum_{k=0}^{m}
\frac{p^{(k)}(y)}{k!}\mu_k^{(n)}(y).
\]
Using $\mu_0^{(n)}=1$, $\mu_1^{(n)}=0$, and the exact formula
for the central moments
\[
\mu_k^{(n)}(y)
=
\sum_{r=\lceil k/2\rceil}^{k-1}
\frac{1}{n^r}
B_{k,k-r}
\bigl(0,\kappa_2(y),\ldots,\kappa_{r+1}(y)\bigr),
\qquad k\geq2,
\]
we may interchange the finite sums and obtain the exact
identity
\begin{equation}
\label{eq:exact-polynomial-expansion}
L_np
=
\sum_{r=0}^{m-1}
\frac{1}{n^r}\mathcal A_r[V]p.
\end{equation}

Apply \eqref{eq:exact-polynomial-expansion} to $p_{x_0}$,
differentiate $\ell$ times with respect to the operator
variable, and then evaluate at $x_0$. By
Lemma~\ref{lem:score-representation},
\begin{equation}
\label{eq:differentiated-polynomial-expansion}
\begin{aligned}
&\sum_{j=0}^{m}
\frac{f^{(j)}(x_0)}{j!}
\mathbb E_{x_0}
\left[
P_{\ell,n}(\eta_{n,x_0},x_0)
\Delta_{n,x_0}^{\,j}
\right]
=
\sum_{r=0}^{m-1}
\frac{1}{n^r}
\bigl(\mathcal A_r[V]p_{x_0}\bigr)^{(\ell)}(x_0).
\end{aligned}
\end{equation}

For $r\geq1$, the coefficient operator has the form
\[
\mathcal A_r[V]g
=
\sum_{k=r+1}^{2r}a_{r,k}g^{(k)},
\]
where the functions $a_{r,k}$ are smooth on $\mathcal M$.
Therefore, $\bigl(\mathcal A_r[V]g\bigr)^{(\ell)}$ involves
derivatives of $g$ of order at most $2r+\ell$.

If $1\leq r\leq s_\ell$, then
\(
2r+\ell\leq2s_\ell+\ell\leq m.
\)
Since $p_{x_0}$ and $f$ have the same derivatives at $x_0$
up to and including order $m$, we obtain
\[
\bigl(\mathcal A_r[V]p_{x_0}\bigr)^{(\ell)}(x_0)
=
\bigl(\mathcal A_r[V]f\bigr)^{(\ell)}(x_0).
\]
For $r=0$, the same equality follows directly from
$\mathcal A_0[V]g=g$ and $\ell\leq m$.

For $s_\ell+1\leq r\leq m-1$, the quantities
\[
\bigl(\mathcal A_r[V]p_{x_0}\bigr)^{(\ell)}(x_0)
\]
are uniformly bounded for $x_0\in J$. Indeed, they are finite
linear combinations of derivatives of the functions
$a_{r,k}$ and derivatives of $p_{x_0}$ at $x_0$; the latter
are equal to $f^{(j)}(x_0)$ for $j\leq m$ and vanish for
$j>m$. Hence,
\[
\sup_{x_0\in J}
\left|
\sum_{r=s_\ell+1}^{m-1}
\frac{1}{n^r}
\bigl(\mathcal A_r[V]p_{x_0}\bigr)^{(\ell)}(x_0)
\right|
\leq
C n^{-s_\ell-1}.
\]
It follows from \eqref{eq:differentiated-polynomial-expansion},
uniformly for $x\in J$, that
\[
\begin{aligned}
&\sum_{j=0}^{m}
\frac{f^{(j)}(x)}{j!}
\mathbb E_x
\left[
P_{\ell,n}(\eta_{n,x},x)\Delta_{n,x}^{\,j}
\right]
=
\sum_{r=0}^{s_\ell}
\frac{1}{n^r}
\bigl(\mathcal A_r[V]f\bigr)^{(\ell)}(x)
+
O\bigl(n^{-s_\ell-1}\bigr).
\end{aligned}
\]
Combining this identity with
\eqref{eq:simultaneous-basic-decomposition} and
\eqref{eq:full-score-remainder}, we obtain
\eqref{eq:quantitative-simultaneous}.
\end{pf}

\begin{corollary}[Simultaneous Expansion of Arbitrary Order]
\label{cor:all-order-simultaneous}
Let $N,\ell\in\mathbb N_0$, and let
$f\colon\mathbb R\to\mathbb R$ be a Borel measurable
function of at most polynomial growth such that
\[
f|_{\mathcal M}\in C^{2N+\ell}(\mathcal M).
\]
Then, for every compact set $J\subset\mathcal M$,
\begin{equation}
\label{eq:all-order-simultaneous}
\sup_{x\in J}
\left|
(L_nf)^{(\ell)}(x)
-
\sum_{r=0}^{N}
\frac{1}{n^r}
\bigl(\mathcal A_r[V]f\bigr)^{(\ell)}(x)
\right|
=
o(n^{-N}).
\end{equation}
\end{corollary}

\begin{pf}
Apply Theorem~\ref{thm:simultaneous-expansion} with
$m=2N+\ell$. Then $s_\ell=N$, and the uniform continuity of
$f^{(2N+\ell)}$ on $J_\rho$ implies that
\[
\omega_{2N+\ell,J_\rho}(f;n^{-1/2})\longrightarrow0.
\]
The second term on the right-hand side of
\eqref{eq:quantitative-simultaneous} is
$O(n^{-N-1})$, which proves
\eqref{eq:all-order-simultaneous}.
\end{pf}

For $\ell=0$, Corollary~\ref{cor:all-order-simultaneous}
reduces to the expansion of arbitrary order established
earlier.

\begin{rmk}[Local Nature of the Theorem]
\label{rem:simultaneous-locality}
The smoothness of $f$ is used only on the compact
neighborhood $J_\rho\Subset\mathcal M$; outside this
neighborhood, the proof requires only Borel measurability and
polynomial growth. Compact containment also guarantees the
uniform boundedness of $V^{-1}$ and the derivatives of $V$
that occur in the score polynomials and coefficient
operators. Thus, no smooth extension of $f$ beyond the mean
domain is required, no assertion is made at the boundary
points of $\mathcal M$, and the constants may deteriorate as
$J$ approaches the boundary.
\end{rmk}

\section{Consequences of the Simultaneous Expansion}
\label{sec:simultaneous-consequences}

\subsection{Order of the Differentiated Coefficient Operators}
\label{subsec:differentiated-coefficient-order}

Corollary~\ref{cor:all-order-simultaneous} shows that the
expansion of the $\ell$th derivative through the term of
order $n^{-N}$ holds under the smoothness assumption
\[
f|_{\mathcal M}\in C^{2N+\ell}(\mathcal M).
\]
We now show that the number $2N+\ell$ is determined by the
structure of the coefficient operators and is not merely a
consequence of the proof method.

\begin{proposition}[Principal Part of the Differentiated
Coefficient Operators]
\label{prop:principal-part-differentiated-coefficients}
Let $N,\ell\in\mathbb N_0$.

For $N\geq1$, there exist functions
\[
b_{N,\ell,j}\in C^\infty(\mathcal M),
\qquad
N+1\leq j\leq2N+\ell-1,
\]
such that, for every function
$f\in C^{2N+\ell}(\mathcal M)$,
\begin{equation}
\label{eq:principal-part-differentiated-coefficient}
\begin{aligned}
\bigl(\mathcal A_N[V]f\bigr)^{(\ell)}(x)
={}&
\frac{V(x)^N}{2^N N!}
f^{(2N+\ell)}(x)
+
\sum_{j=N+1}^{2N+\ell-1}
b_{N,\ell,j}(x)f^{(j)}(x),
\qquad x\in\mathcal M.
\end{aligned}
\end{equation}
Here and below, an empty sum is understood to be zero.

For $N=0$,
\[
\bigl(\mathcal A_0[V]f\bigr)^{(\ell)}
=
f^{(\ell)}.
\]
Consequently, the differential operator
\[
f\longmapsto
\bigl(\mathcal A_N[V]f\bigr)^{(\ell)}
\]
has exact order $2N+\ell$.
\end{proposition}

\begin{pf}
For $N=0$, the assertion follows from
$\mathcal A_0[V]f=f$. Let $N\geq1$. Differentiating
\eqref{eq:coefficient-operator-principal-part} $\ell$ times
and applying the Leibniz rule yields
\eqref{eq:principal-part-differentiated-coefficient}.
Since $V>0$ on $\mathcal M$, the leading coefficient
\(
\frac{V^N}{2^N N!}
\)
does not vanish. Therefore, the resulting differential
operator has exact order $2N+\ell$.
\end{pf}

\begin{corollary}[Minimal Classical Smoothness]
\label{cor:minimal-classical-smoothness}
Let $N,\ell\in\mathbb N_0$. In the class of Borel
measurable functions
\(
f\colon\mathbb R\to\mathbb R
\)
of at most polynomial growth, the condition
\[
f|_{\mathcal M}
\in
C^{2N+\ell}(\mathcal M)
\]
is sufficient for the expansion of $(L_nf)^{(\ell)}$
through the term of order $n^{-N}$, uniformly on compact
subsets. Moreover, $2N+\ell$ is the smallest integer order
of smoothness guaranteeing that every coefficient
\[
\bigl(\mathcal A_r[V]f\bigr)^{(\ell)},
\qquad
0\leq r\leq N,
\]
is classically defined for every function in this class.
\end{corollary}

\begin{pf}
Sufficiency follows directly from
Corollary~\ref{cor:all-order-simultaneous}. By
Proposition~\ref{prop:principal-part-differentiated-coefficients},
the last coefficient contains the nonzero term
\(
\frac{V^N}{2^N N!}f^{(2N+\ell)}.
\)
Thus, if $m<2N+\ell$, this coefficient need not be
classically defined for a function in $C^m(\mathcal M)$.
\end{pf}

For comparison, Abel's general result on the termwise
differentiation of full asymptotic expansions for operators
of exponential type~\cite{Abel2023} requires smoothness of
order $2(N+\ell)$ for the expansion of the $\ell$th
derivative through the term of order $n^{-N}$. In the
present proof, the score-polynomial representation
\eqref{eq:score-representation}, together with the exact
polynomial identity
\eqref{eq:exact-polynomial-expansion}, makes it possible to
identify the differentiated coefficients without
differentiating the asymptotic remainder. This yields the
sufficient classical smoothness order $2N+\ell$, saving
$\ell$ derivatives relative to Abel's condition.

The sharpness of this threshold for the asymptotic expansion
itself, rather than merely for the classical formula for its
coefficients, will be established in
Proposition~\ref{prop:gaussian-smoothness-sharpness}.

The smoothness requirement naturally splits into two parts:
\[
\underbrace{2N}_{\text{expansion through }n^{-N}}
+
\underbrace{\ell}_{\ell\text{ simultaneous differentiations}}.
\]
Thus, each successive asymptotic order requires two
additional derivatives, whereas each additional
differentiation in the simultaneous approximation requires
one additional derivative. The minimality assertion is
understood within the classical uniform scale $C^m$; it does
not exclude weaker local assumptions for particular
functions or expansions formulated at an individual point.

If two additional derivatives are available, then applying
Corollary~\ref{cor:all-order-simultaneous} with $N+1$ gives
\[
(L_nf)^{(\ell)}
-
\sum_{r=0}^{N}
\frac{1}{n^r}
\bigl(\mathcal A_r[V]f\bigr)^{(\ell)}
=
\frac{1}{n^{N+1}}
\bigl(\mathcal A_{N+1}[V]f\bigr)^{(\ell)}
+
o(n^{-N-1})
\]
locally uniformly on $\mathcal M$. In particular, if
$f|_{\mathcal M}\in C^{2N+\ell+2}(\mathcal M)$, then, after
truncation following the term of order $n^{-N}$, the
remainder is of order $O(n^{-N-1})$.

\subsection{Explicit Lower-Order Formulas}
\label{subsec:explicit-simultaneous-formulas}

The first nontrivial coefficient operator is
\(
\mathcal A_1[V]f
=
\frac{V}{2}f''.
\)
Therefore, the general simultaneous expansion yields
particularly simple explicit formulas for the first two
derivatives of $L_nf$.

\begin{corollary}[First-Order Simultaneous Formulas]
\label{cor:first-order-simultaneous-formulas}
Let $f\colon\mathbb R\to\mathbb R$ be a Borel measurable
function of at most polynomial growth.

\begin{enumerate}
\item
If
\(
f|_{\mathcal M}\in C^3(\mathcal M),
\)
then
\[
(L_nf)'(x)
=
f'(x)
+
\frac{1}{2n}
\bigl(Vf''\bigr)'(x)
+
o(n^{-1})
\]
locally uniformly on $\mathcal M$. Equivalently,
\[
(L_nf)'(x)
=
f'(x)
+
\frac{1}{2n}
\left(
V'(x)f''(x)+V(x)f'''(x)
\right)
+
o(n^{-1}).
\]

\item
If
\(
f|_{\mathcal M}\in C^4(\mathcal M),
\)
then
\[
(L_nf)''(x)
=
f''(x)
+
\frac{1}{2n}
\bigl(Vf''\bigr)''(x)
+
o(n^{-1})
\]
locally uniformly on $\mathcal M$. Equivalently,
\[
\begin{aligned}
(L_nf)''(x)
={}&
f''(x)
+
\frac{1}{2n}
\Bigl(
V''(x)f''(x)
+
2V'(x)f'''(x)
+
V(x)f^{(4)}(x)
\Bigr)
+
o(n^{-1}).
\end{aligned}
\]
\end{enumerate}
\end{corollary}

\begin{pf}
Apply Corollary~\ref{cor:all-order-simultaneous} with
$(N,\ell)=(1,1)$ and $(N,\ell)=(1,2)$, respectively, and
use
\(
\mathcal A_0[V]f=f,
\) 
\(
\mathcal A_1[V]f=\frac{V}{2}f''.
\)
The explicit formulas follow from the Leibniz rule.
\end{pf}

\section{Sharpness Results}
\label{sec:sharpness}

Throughout this section, we consider the Gaussian natural
exponential family for which
\(
\mathcal M=\mathbb R,
\) 
\(
V=\frac12.
\)
The corresponding sample-mean operators are the
Gauss--Weierstrass operators
\[
W_nf(x)
=
\sqrt{\frac{n}{\pi}}
\int_{\mathbb R}
f(t)e^{-n(t-x)^2}\,dt.
\]

Since $\kappa_2=1/2$ and $\kappa_j=0$ for $j\geq3$, the
multi-index representation
\eqref{eq:coefficient-operator-multi-index}, together with
the identity $\mathcal A_0[V]f=f$, gives
\begin{equation}
\label{eq:gaussian-coefficients-sharpness}
\mathcal A_r\left[\frac12\right]f
=
\frac{1}{4^r r!}f^{(2r)},
\qquad r\in\mathbb N_0.
\end{equation}

Fix a function $\chi\in C_c^\infty(\mathbb R)$ such that
\[
0\leq\chi\leq1,
\qquad
\chi=1\quad\text{on }[-1,1],
\qquad
\chi=0\quad\text{outside }[-2,2],
\]
and, for $\beta>0$, set
\[
h_\beta(t)
=
\chi(t)(t_+)^\beta,
\qquad
t_+=\max\{t,0\}.
\]

The following elementary scaling lemma will be used in both
sharpness proofs.

\begin{lemma}[Gaussian Scaling]
\label{lem:gaussian-scaling}
Let $\ell\in\mathbb N_0$ and $\beta>\ell$. Then
\begin{equation}
\label{eq:gaussian-scaling}
(W_nh_\beta)^{(\ell)}(0)
\sim
\frac{\Gamma(\beta+1)}
     {\Gamma(\beta-\ell+1)}
\frac{
\Gamma\left(\frac{\beta-\ell+1}{2}\right)
}{2\sqrt{\pi}}\,
n^{-(\beta-\ell)/2},
\qquad n\to\infty.
\end{equation}
\end{lemma}

\begin{pf}
Since $h_\beta$ is compactly supported and belongs to
$C^\ell(\mathbb R)$, differentiation under the integral
sign followed by integration by parts gives
\[
(W_nh_\beta)^{(\ell)}(0)
=
\sqrt{\frac{n}{\pi}}
\int_{\mathbb R}
h_\beta^{(\ell)}(t)e^{-nt^2}\,dt.
\]
For $0<t<1$, we have
\[
h_\beta^{(\ell)}(t)
=
\frac{\Gamma(\beta+1)}
     {\Gamma(\beta-\ell+1)}
t^{\beta-\ell}.
\]
On the remaining part of its support, namely on $[1,2]$,
the function $h_\beta^{(\ell)}$ is bounded. Therefore,
\[
\begin{aligned}
(W_nh_\beta)^{(\ell)}(0)
={}&
\frac{\Gamma(\beta+1)}
     {\Gamma(\beta-\ell+1)}
\sqrt{\frac{n}{\pi}}
\int_0^1
t^{\beta-\ell}e^{-nt^2}\,dt
+
O\bigl(\sqrt n\,e^{-n}\bigr)
\\
={}&
\frac{\Gamma(\beta+1)}
     {\Gamma(\beta-\ell+1)}
\frac{n^{-(\beta-\ell)/2}}{\sqrt{\pi}}
\int_0^{\sqrt n}
u^{\beta-\ell}e^{-u^2}\,du
+
O\bigl(\sqrt n\,e^{-n}\bigr).
\end{aligned}
\]
Finally,
\[
\int_0^\infty
u^{\beta-\ell}e^{-u^2}\,du
=
\frac12
\Gamma\left(\frac{\beta-\ell+1}{2}\right),
\]
which proves \eqref{eq:gaussian-scaling}.
\end{pf}

\subsection{Sharpness of the Smoothness Threshold}
\label{subsec:smoothness-sharpness}

We first show that the smoothness threshold
$C^{2N+\ell}$ in
Corollary~\ref{cor:all-order-simultaneous} cannot, in
general, be lowered within the scale of integer-order
smoothness classes.

\begin{proposition}[Sharpness in the Gaussian Case]
\label{prop:gaussian-smoothness-sharpness}
Let $N\in\mathbb N$, $\ell\in\mathbb N_0$, and
$0<\alpha<1$. Define
\(
f_\alpha
=
h_{2N+\ell-1+\alpha}.
\)
Then
\[
f_\alpha
\in
C^{2N+\ell-1}(\mathbb R)
\setminus
C^{2N+\ell}(\mathbb R),
\]
and
\begin{equation}
\label{eq:gaussian-smoothness-sharpness}
(W_nf_\alpha)^{(\ell)}(0)
\sim
C_{N,\ell,\alpha}
n^{-N+(1-\alpha)/2},
\qquad n\to\infty,
\end{equation}
where
\[
C_{N,\ell,\alpha}
=
\frac{\Gamma(2N+\ell+\alpha)}
     {\Gamma(2N+\alpha)}
\frac{\Gamma(N+\alpha/2)}
     {2\sqrt{\pi}}
>0.
\]
Moreover,
\[
\bigl(\mathcal A_r[1/2]f_\alpha\bigr)^{(\ell)}(0)
=
0,
\qquad
0\leq r\leq N-1.
\]
Consequently,
\begin{equation}
\label{eq:gaussian-smoothness-divergence}
\begin{aligned}
n^N
\left|
(W_nf_\alpha)^{(\ell)}(0)
-
\sum_{r=0}^{N-1}
\frac{1}{n^r}
\bigl(\mathcal A_r[1/2]f_\alpha\bigr)^{(\ell)}(0)
\right|
\longrightarrow\infty.
\end{aligned}
\end{equation}
Thus, $C^{2N+\ell-1}$ smoothness is insufficient for an
expansion with remainder $o(n^{-N})$.
\end{proposition}

\begin{pf}
Set
\(
k=2N+\ell-1.
\)
Since
\(
f_\alpha(t)
=
\chi(t)(t_+)^{k+\alpha}
\)
and $0<\alpha<1$, we have
\(
f_\alpha
\in
C^k(\mathbb R)\setminus C^{k+1}(\mathbb R).
\)
Moreover,
\(
f_\alpha^{(j)}(0)=0,
\) 
\(
0\leq j\leq k.
\)

Apply Lemma~\ref{lem:gaussian-scaling} with
\(
\beta
=
2N+\ell-1+\alpha.
\)
Since
\(
\beta-\ell
=
2N-1+\alpha,
\)
formula \eqref{eq:gaussian-scaling} gives
\[
(W_nf_\alpha)^{(\ell)}(0)
\sim
C_{N,\ell,\alpha}
n^{-N+(1-\alpha)/2},
\]
which proves
\eqref{eq:gaussian-smoothness-sharpness}.

By \eqref{eq:gaussian-coefficients-sharpness},
\[
\bigl(\mathcal A_r[1/2]f_\alpha\bigr)^{(\ell)}(0)
=
\frac{1}{4^r r!}
f_\alpha^{(2r+\ell)}(0).
\]
For $0\leq r\leq N-1$,
\[
2r+\ell
\leq
2N+\ell-2
\leq k,
\]
and hence all these coefficients vanish. Therefore, the
expression inside the absolute value in
\eqref{eq:gaussian-smoothness-divergence} is equal to
$(W_nf_\alpha)^{(\ell)}(0)$. It follows that
\[
\begin{aligned}
n^N
\left|(W_nf_\alpha)^{(\ell)}(0)\right|
&\sim
C_{N,\ell,\alpha}
n^{(1-\alpha)/2}
\longrightarrow\infty.
\end{aligned}
\]
This proves
\eqref{eq:gaussian-smoothness-divergence}.
\end{pf}

The function $f_\alpha$ is compactly supported and therefore
satisfies the global assumptions of Borel measurability and
polynomial growth in the simultaneous expansion theorem.
Thus, the obstruction is purely local and is caused precisely
by the absence of a derivative of order $2N+\ell$.

\subsection{Sharpness of the Quantitative Remainder}
\label{subsec:quantitative-remainder-sharpness}

The preceding counterexample concerns the integer-order
smoothness threshold. We now show, without any parity
restriction on $m-\ell$, that the modulus-of-continuity term
in Theorem~\ref{thm:simultaneous-expansion} has the correct
asymptotic order.

\begin{proposition}[Sharpness of the Quantitative Remainder]
\label{prop:quantitative-remainder-sharpness}
Let $m,\ell\in\mathbb N_0$, $m\geq\ell$, and
$0<\alpha<1$. Set
\(
s_\ell
=
\left\lfloor\frac{m-\ell}{2}\right\rfloor
\)
and define
\(
g_\alpha
=
h_{m+\alpha}.
\)
Then
\[
g_\alpha
\in
C^m(\mathbb R)\setminus C^{m+1}(\mathbb R).
\]

Let $J=\{0\}$ and $\rho=1/2$, so that
\(
J_\rho=[-1/2,1/2].
\)
For every $0<\delta\leq1/2$,
\begin{equation}
\label{eq:sharp-modulus}
\omega_{m,J_\rho}(g_\alpha;\delta)
=
c_{m,\alpha}\delta^\alpha,
\qquad
c_{m,\alpha}
=
\frac{\Gamma(m+\alpha+1)}
     {\Gamma(\alpha+1)}.
\end{equation}
Moreover,
\begin{equation}
\label{eq:sharp-quantitative-asymptotic}
(W_ng_\alpha)^{(\ell)}(0)
\sim
d_{m,\ell,\alpha}
n^{-(m-\ell+\alpha)/2},
\qquad n\to\infty,
\end{equation}
where
\[
d_{m,\ell,\alpha}
=
\frac{\Gamma(m+\alpha+1)}
     {\Gamma(m-\ell+\alpha+1)}
\frac{
\Gamma\left(\frac{m-\ell+\alpha+1}{2}\right)
}{2\sqrt{\pi}}
>0.
\]
Furthermore,
\[
\bigl(\mathcal A_r[1/2]g_\alpha\bigr)^{(\ell)}(0)
=
0,
\qquad
0\leq r\leq s_\ell,
\]
and hence
\begin{equation}
\label{eq:sharp-quantitative-comparison}
\begin{aligned}
&
\left|
(W_ng_\alpha)^{(\ell)}(0)
-
\sum_{r=0}^{s_\ell}
\frac{1}{n^r}
\bigl(\mathcal A_r[1/2]g_\alpha\bigr)^{(\ell)}(0)
\right|
\sim
\frac{d_{m,\ell,\alpha}}{c_{m,\alpha}}\,
n^{-(m-\ell)/2}
\omega_{m,J_\rho}
\bigl(g_\alpha;n^{-1/2}\bigr).
\end{aligned}
\end{equation}
Thus, for every $m\geq\ell$, the
modulus-of-continuity term in
\eqref{eq:quantitative-simultaneous} is sharp in order.
\end{proposition}

\begin{pf}
Since
\(
g_\alpha(t)
=
\chi(t)(t_+)^{m+\alpha}
\)
and $0<\alpha<1$, we have
\[
g_\alpha
\in
C^m(\mathbb R)\setminus C^{m+1}(\mathbb R),
\qquad
g_\alpha^{(j)}(0)=0,
\quad
0\leq j\leq m.
\]

On $J_\rho=[-1/2,1/2]$, the cutoff function is identically
equal to one, and therefore
\(
g_\alpha^{(m)}(t)
=
c_{m,\alpha}(t_+)^\alpha.
\)
The elementary inequality
\(
|u_+^\alpha-v_+^\alpha|
\leq
|u-v|^\alpha,
\)
\(
u,v\in\mathbb R,
\)
shows that
\[
\omega_{m,J_\rho}(g_\alpha;\delta)
\leq
c_{m,\alpha}\delta^\alpha.
\]
Choosing $u=\delta$ and $v=0$ gives the reverse inequality.
Hence, \eqref{eq:sharp-modulus} holds.

Apply Lemma~\ref{lem:gaussian-scaling} with
\(
\beta
=
m+\alpha.
\)
Since
\(
\beta-\ell
=
m-\ell+\alpha,
\)
formula \eqref{eq:gaussian-scaling} gives
\[
(W_ng_\alpha)^{(\ell)}(0)
\sim
d_{m,\ell,\alpha}
n^{-(m-\ell+\alpha)/2},
\]
which proves
\eqref{eq:sharp-quantitative-asymptotic}.

By \eqref{eq:gaussian-coefficients-sharpness},
\[
\bigl(\mathcal A_r[1/2]g_\alpha\bigr)^{(\ell)}(0)
=
\frac{1}{4^r r!}
g_\alpha^{(2r+\ell)}(0).
\]
For $0\leq r\leq s_\ell$,
\[
2r+\ell
\leq
2s_\ell+\ell
\leq
m.
\]
Consequently,
\(
g_\alpha^{(2r+\ell)}(0)=0,
\)
\(
0\leq r\leq s_\ell,
\)
and all the coefficients in
\eqref{eq:sharp-quantitative-comparison} vanish. Therefore,
the expression inside the absolute value in
\eqref{eq:sharp-quantitative-comparison} is equal to
\[
\left|(W_ng_\alpha)^{(\ell)}(0)\right|.
\]

Finally, \eqref{eq:sharp-modulus} gives
\[
\begin{aligned}
&
n^{-(m-\ell)/2}
\omega_{m,J_\rho}
\bigl(g_\alpha;n^{-1/2}\bigr)
=
c_{m,\alpha}
n^{-(m-\ell+\alpha)/2}.
\end{aligned}
\]
Combining this identity with
\eqref{eq:sharp-quantitative-asymptotic}, we obtain
\eqref{eq:sharp-quantitative-comparison}.
\end{pf}

Since $0<\alpha<1$,
\[
\frac{m-\ell+\alpha}{2}
<
\left\lfloor\frac{m-\ell}{2}\right\rfloor+1
=
s_\ell+1.
\]
It follows that
\[
n^{-s_\ell-1}
=
o\left(
n^{-(m-\ell+\alpha)/2}
\right).
\]
Thus, the second term in
\eqref{eq:quantitative-simultaneous} is negligible for this
example, whereas the modulus-of-continuity term has exactly
the same order as the approximation error. Consequently, for every $m\geq\ell$, the order of the
modulus-of-continuity term in
\eqref{eq:quantitative-simultaneous} is optimal over the class
of functions covered by
Theorem~\ref{thm:simultaneous-expansion}. In general, this term
cannot be replaced by
\[
o\left(
n^{-(m-\ell)/2}
\omega_{m,J_\rho}\bigl(f;n^{-1/2}\bigr)
\right).
\]
\section{Applications to Classical Positive Operators}
\label{sec:classical-operators}

The sample-mean construction encompasses a number of
classical positive approximation operators. The five
families considered below are listed in
Table~\ref{tab:classical-operators}. For general information
on positive linear approximation operators and
exponential-type constructions, see
\cite{Paltanea2004} and \cite{IsmailMay1978}, respectively.

\begin{table}[htbp]
\centering
\caption{Natural exponential families and the corresponding
sample-mean operators}\label{tab:classical-operators}
\begin{tabular*}{.9\linewidth}{@{\extracolsep{\fill}}llll@{}}
\toprule
Distribution of $\xi_x$
&
Mean domain $\mathcal M$
&
$V(x)$
&
Operator
\\
\midrule
Bernoulli
&
$(0,1)$
&
$x(1-x)$
&
Bernstein $B_n$
\\
Poisson
&
$(0,\infty)$
&
$x$
&
Sz\'asz--Mirakjan $S_n$
\\
Geometric
&
$(0,\infty)$
&
$x(1+x)$
&
Baskakov $\mathcal B_n$
\\
Exponential
&
$(0,\infty)$
&
$x^2$
&
Post--Widder $F_n$
\\
Normal $N(x,1/2)$
&
$\mathbb R$
&
$1/2$
&
Gauss--Weierstrass $W_n$
\\
\bottomrule
\end{tabular*}
\end{table}

Throughout this section, we assume that
$f\colon\mathbb R\to\mathbb R$ is Borel measurable, is of at
most polynomial growth, and satisfies
\(
f|_{\mathcal M}\in C^4(\mathcal M)
\)
for the corresponding mean domain. Applying
Corollary~\ref{cor:ordinary-all-order-expansion} with $N=2$
and using \eqref{eq:first-coefficient-operators}, we obtain
the common second-order formula
\begin{equation}
\label{eq:classical-second-order-template}
\begin{aligned}
L_nf(x)
={}&
f(x)
+
\frac{V(x)}{2n}f''(x)
+
\frac{1}{n^2}
\left[
\frac{V(x)V'(x)}{6}f'''(x)
+
\frac{V(x)^2}{8}f^{(4)}(x)
\right]
+
o\left(n^{-2}\right),
\end{aligned}
\end{equation}
locally uniformly on $\mathcal M$. Thus, in each example it
suffices to specify the operator and substitute the
corresponding variance function.

\subsection{Discrete Operators}
\label{subsec:discrete-operators}

\paragraph{Bernstein operators.}
The classical Bernstein operators~\cite{Bernstein1912},
corresponding to the Bernoulli family, are given by
\[
B_nf(x)
=
\sum_{k=0}^{n}
f\left(\frac{k}{n}\right)
\binom{n}{k}
x^k(1-x)^{n-k}.
\]
Substituting $V(x)=x(1-x)$ into
\eqref{eq:classical-second-order-template}, we obtain
\[
\begin{aligned}
B_nf(x)
={}&
f(x)
+
\frac{x(1-x)}{2n}f''(x)\\
&
+
\frac{1}{n^2}
\left[
\frac{x(1-x)(1-2x)}{6}f'''(x)
+
\frac{x^2(1-x)^2}{8}f^{(4)}(x)
\right]
+
o\left(n^{-2}\right).
\end{aligned}
\]

This is the truncation of the full Bernstein expansion after
the second-order term; see the work of
Tachev~\cite{Tachev2012}.

\paragraph{Sz\'asz--Mirakjan operators.}
The Poisson family gives rise to the Sz\'asz--Mirakjan
operators~\cite{Mirakjan1941,Szasz1950}:
\[
S_nf(x)
=
e^{-nx}
\sum_{k=0}^{\infty}
\frac{(nx)^k}{k!}
f\left(\frac{k}{n}\right).
\]
For $V(x)=x$, formula
\eqref{eq:classical-second-order-template} becomes
\[
\begin{aligned}
S_nf(x)
={}&
f(x)
+
\frac{x}{2n}f''(x)
+
\frac{1}{n^2}
\left[
\frac{x}{6}f'''(x)
+
\frac{x^2}{8}f^{(4)}(x)
\right]
+
o\left(n^{-2}\right).
\end{aligned}
\]

\paragraph{Baskakov operators.}
For the geometric family, the
Baskakov operators~\cite{Baskakov1957} are given by
\[
\mathcal B_nf(x)
=
\sum_{k=0}^{\infty}
\binom{n+k-1}{k}
\frac{x^k}{(1+x)^{n+k}}
f\left(\frac{k}{n}\right).
\]
Since $V(x)=x(1+x)$,
\[
\begin{aligned}
\mathcal B_nf(x)
={}&
f(x)
+
\frac{x(1+x)}{2n}f''(x)\\
&
+
\frac{1}{n^2}
\left[
\frac{x(1+x)(1+2x)}{6}f'''(x)
+
\frac{x^2(1+x)^2}{8}f^{(4)}(x)
\right]
+
o\left(n^{-2}\right).
\end{aligned}
\]

\subsection{Integral Operators}
\label{subsec:integral-operators}

\paragraph{Post--Widder operators.}
The exponential family gives rise to the Post--Widder
operators~\cite{Post1930,Widder1934}:
\[
F_nf(x)
=
\frac{n^n}{\Gamma(n)x^n}
\int_0^\infty
f(t)t^{n-1}e^{-nt/x}\,dt.
\]
For $V(x)=x^2$, we obtain
\[
\begin{aligned}
F_nf(x)
={}&
f(x)
+
\frac{x^2}{2n}f''(x)
+
\frac{1}{n^2}
\left[
\frac{x^3}{3}f'''(x)
+
\frac{x^4}{8}f^{(4)}(x)
\right]
+
o\left(n^{-2}\right).
\end{aligned}
\]

A full pointwise asymptotic expansion for the classical
Post--Widder operators was also obtained as a special case
of the result of Abel, Agratini, and P\u{a}lt\u{a}nea for
semi-exponential Post--Widder
operators~\cite{AbelAgratiniPaltanea2026}.

\paragraph{Gauss--Weierstrass operators.}
For the normal family with variance $1/2$,
\[
W_nf(x)
=
\sqrt{\frac{n}{\pi}}
\int_{\mathbb R}
f(t)e^{-n(t-x)^2}\,dt.
\]
Using \eqref{eq:gaussian-coefficients-sharpness}, we obtain
\[
W_nf(x)
=
f(x)
+
\frac{1}{4n}f''(x)
+
\frac{1}{32n^2}f^{(4)}(x)
+
o\left(n^{-2}\right).
\]

Thus, all five examples are realizations of the same
sample-mean construction. Their coefficients of arbitrary
order are obtained from \eqref{eq:coefficient-operator}
using the variance function given in
Table~\ref{tab:classical-operators}.

\subsection{Examples of Simultaneous Expansions}
\label{subsec:simultaneous-examples}

The simultaneous expansion also yields explicit asymptotic
formulas for derivatives of classical operators. We present
two representative cases.

\paragraph{First derivative of the Bernstein operators.}
Let $f\colon\mathbb R\to\mathbb R$ be a Borel measurable
function of at most polynomial growth such that
\[
f|_{(0,1)}\in C^3(0,1).
\]
Applying Corollary~\ref{cor:all-order-simultaneous} with
$N=1$, $\ell=1$, and $V(x)=x(1-x)$, we obtain
\[
\begin{aligned}
(B_nf)'(x)
={}&
f'(x)
+
\frac{1}{2n}
\bigl(x(1-x)f''(x)\bigr)'
+
o(n^{-1})\\
&
={}
f'(x)
+
\frac{1}{2n}
\left[
(1-2x)f''(x)
+
x(1-x)f'''(x)
\right]
+
o(n^{-1}),
\end{aligned}
\]
locally uniformly for $x\in(0,1)$.

\paragraph{Arbitrary-order derivatives of the
Gauss--Weierstrass operators.}
Let $N,\ell\in\mathbb N_0$, and let
$f\colon\mathbb R\to\mathbb R$ be a Borel measurable
function of at most polynomial growth such that
\[
f\in C^{2N+\ell}(\mathbb R).
\]
From \eqref{eq:gaussian-coefficients-sharpness} and
Corollary~\ref{cor:all-order-simultaneous}, we obtain
\begin{equation}
\label{eq:gauss-weierstrass-simultaneous-expansion}
(W_nf)^{(\ell)}(x)
=
\sum_{r=0}^{N}
\frac{1}{4^r r!\,n^r}
f^{(2r+\ell)}(x)
+
o(n^{-N})
\end{equation}
locally uniformly for $x\in\mathbb R$.

Full asymptotic expansions for the Gauss--Weierstrass
operators, including simultaneous approximation, were also
obtained by Abel, Agratini, and
P\u{a}lt\u{a}nea~\cite{AbelAgratiniPaltanea2018}.

The Gaussian family is also the model used in
Section~\ref{sec:sharpness}. In particular,
Propositions~\ref{prop:gaussian-smoothness-sharpness} and
\ref{prop:quantitative-remainder-sharpness} show that the
smoothness threshold and the order of the
modulus-of-continuity term underlying
\eqref{eq:gauss-weierstrass-simultaneous-expansion} are, in
general, sharp.





\subsection*{Funding}

This research received no external funding.

\subsection*{Data availability}

No data were used for the research described in this article.

\subsection*{Declaration of Competing Interest}

The author declares no competing interests.


\begin{thebibliography}{99}

\bibitem{Abel2023}
U.~Abel, \emph{Simultaneous approximation by operators of exponential type},
arXiv:2309.09099 [math.CA], 2023.

\bibitem{AbelAgratiniPaltanea2018}
U.~Abel, O.~Agratini and R.~P{\u{a}}lt{\u{a}}nea,
\emph{A complete asymptotic expansion for the quasi-interpolants of
Gau{\ss}-Weierstra{\ss} operators},
Mediterr. J. Math. \textbf{15} (2018), no.~4, Paper No.~156, 10~pp.
DOI \texttt{10.1007/s00009-018-1195-8}.

\bibitem{AbelAgratiniPaltanea2026}
U.~Abel, O.~Agratini and R.~P{\u{a}}lt{\u{a}}nea,
\emph{A complete asymptotic expansion for the semi-exponential
Post-Widder operators},
Adv. Oper. Theory \textbf{11} (2026), no.~2, Paper No.~24, 15~pp.
DOI \texttt{10.1007/s43036-026-00499-5}.

\bibitem{AbelGupta2021}
U.~Abel and V.~Gupta,
\emph{A complete asymptotic expansion for operators of exponential type
with $p(x)=x(1+x)^2$},
Positivity \textbf{25} (2021), no.~3, 1013--1025.
DOI \texttt{10.1007/s11117-020-00802-5}.

\bibitem{AbelGuptaKushnirevych2021}
U.~Abel, V.~Gupta and V.~Kushnirevych,
\emph{Asymptotic expansions for certain exponential-type operators
connected with $2x^{3/2}$},
Math. Sci. \textbf{15} (2021), no.~3, 311--315.
DOI \texttt{10.1007/s40096-021-00382-9}.

\bibitem{AbelKushnirevych2025}
U.~Abel and V.~Kushnirevych,
\emph{Complete asymptotic expansion for exponential-type operators
related to $x^{4/3}$},
Math. Pannon. (N.~S.) \textbf{31} (2025), no.~1, 68--76.
DOI \texttt{10.1556/314.2025.00008}.

\bibitem{BarndorffNielsen1978}
O.~Barndorff-Nielsen,
\emph{Information and exponential families in statistical theory},
Wiley Series in Probability and Statistics,
John Wiley \& Sons, Ltd., Chichester, 2014.
Reprint of the 1978 original.
DOI \texttt{10.1002/9781118857281}.

\bibitem{Baskakov1957}
V.~A. Baskakov,
\emph{An instance of a sequence of linear positive operators in the
space of continuous functions},
Dokl. Akad. Nauk SSSR (N.S.) \textbf{113} (1957), 249--251.

\bibitem{Bernstein1912}
S.~Bernstein,
\emph{D\'emonstration du th\'eor\`eme de Weierstrass fond\'ee sur le calcul
des probabilit\'es},
Communications of the Kharkov Mathematical Society \textbf{13} (1912), 1--2.

\bibitem{Comtet1974}
L.~Comtet,
\emph{Advanced combinatorics. The art of finite and infinite expansions},
enlarged ed., D.~Reidel Publishing Co., Dordrecht, 1974.

\bibitem{Floater2005}
M.~S. Floater,
\emph{On the convergence of derivatives of Bernstein approximation},
J. Approx. Theory \textbf{134} (2005), no.~1, 130--135.
DOI \texttt{10.1016/j.jat.2004.02.009}.

\bibitem{GonskaPaltanea2010}
H.~Gonska and R.~P{\u{a}}lt{\u{a}}nea,
\emph{General Voronovskaya and asymptotic theorems in simultaneous
approximation},
Mediterr. J. Math. \textbf{7} (2010), no.~1, 37--49.
DOI \texttt{\nolinkurl{10.1007/s00009-010-0025-4}}.

\bibitem{Holhos2019}
A.~Holho{\c{s}},
\emph{A Voronovskaya-type theorem for the first derivatives of positive
linear operators},
Results Math. \textbf{74} (2019), no.~2, Paper No.~76, 13~pp.
DOI \texttt{10.1007/s00025-019-0992-0}.

\bibitem{Holhos2022Simultaneous}
A.~Holho{\c{s}},
\emph{A Voronovskaya-type theorem in simultaneous approximation},
Period. Math. Hungar. \textbf{85} (2022), no.~2, 280--291.
DOI \texttt{10.1007/s10998-021-00434-6}.

\bibitem{Holhos2022Exponential}
A.~Holho{\c{s}},
\emph{Voronovskaya-type results for positive linear operators of
exponential type and their derivatives},
Bull. Malays. Math. Sci. Soc. \textbf{45} (2022), no.~4, 1839--1861.
DOI \texttt{10.1007/s40840-021-01227-3}.

\bibitem{IsmailMay1978}
M.~E.~H. Ismail and C.~P. May,
\emph{On a family of approximation operators},
J. Math. Anal. Appl. \textbf{63} (1978), no.~2, 446--462.
DOI \texttt{10.1016/0022-247X(78)90090-2}.

\bibitem{Letac1992}
G.~Letac,
\emph{Lectures on natural exponential families and their variance functions},
Monograf\'ias de Matem\'atica [Mathematical Monographs], vol.~50,
Instituto de Matem\'atica Pura e Aplicada (IMPA), Rio de Janeiro, 1992.

\bibitem{May1976}
C.~P. May,
\emph{Saturation and inverse theorems for combinations of a class of
exponential-type operators},
Canadian J. Math. \textbf{28} (1976), no.~6, 1224--1250.
DOI \texttt{10.4153/CJM-1976-123-8}.

\bibitem{Mirakjan1941}
G.~M. Mirakjan,
\emph{Approximation of Continuous Functions with the Aid of Polynomials},
Doklady Akademii Nauk SSSR \textbf{31} (1941), 201--205.

\bibitem{Morris1982}
C.~N. Morris,
\emph{Natural exponential families with quadratic variance functions},
Ann. Statist. \textbf{10} (1982), no.~1, 65--80.

\bibitem{Paltanea2004}
R.~P{\u{a}}lt{\u{a}}nea,
\emph{Approximation theory using positive linear operators},
Birkh\"auser Boston, Inc., Boston, MA, 2004.
DOI \texttt{10.1007/978-1-4612-2058-9}.

\bibitem{Post1930}
E.~L. Post,
\emph{Generalized differentiation},
Trans. Amer. Math. Soc. \textbf{32} (1930), no.~4, 723--781.
DOI \texttt{10.2307/1989348}.

\bibitem{Sikkema1970}
P.~C. Sikkema,
\emph{On some linear positive operators},
Indag. Math. \textbf{32} (1970), 327--337
(Nederl. Akad. Wetensch. Proc. Ser.~A \textbf{73}).

\bibitem{Szasz1950}
O.~Szasz,
\emph{Generalization of S.~Bernstein's polynomials to the infinite interval},
J. Research Nat. Bur. Standards \textbf{45} (1950), 239--245.

\bibitem{Tachev2012}
G.~T. Tachev,
\emph{The complete asymptotic expansion for Bernstein operators},
J. Math. Anal. Appl. \textbf{385} (2012), no.~2, 1179--1183.
DOI \texttt{10.1016/j.jmaa.2011.07.042}.

\bibitem{Voronovskaya1932}
E.~V. Voronovskaja,
\emph{D\'etermination de la forme asymptotique d'approximation des
fonctions par les polyn\^omes de M.~Bernstein},
C. R. Acad. Sci. URSS \textbf{4} (1932), 79--85.

\bibitem{Widder1934}
D.~V. Widder,
\emph{The inversion of the Laplace integral and the related moment problem},
Trans. Amer. Math. Soc. \textbf{36} (1934), no.~1, 107--200.
DOI \texttt{10.2307/1989710}.

\bibitem{Xiang2022}
J.~X. Xiang,
\emph{On convergence of derivatives of Voronovskaja-type expansion for
Bernstein type operators},
J. Approx. Theory \textbf{273} (2022), Paper No.~105668, 19~pp.
DOI \texttt{10.1016/j.jat.2021.105668}.

\end{thebibliography}
\end{document}